\documentclass[letterpaper, 11pt]{article}

\usepackage{amsmath, amssymb}
\usepackage[T1]{fontenc}
\usepackage[left=2.5cm, right=2.5cm, top=2.5cm, bottom = 2.5cm]{geometry}
\usepackage[utf8]{inputenc}
\usepackage{makecell}
\usepackage{multirow}
\usepackage{placeins}
\usepackage[center]{subfigure}
\usepackage{tikz}
\usepackage{doi}
\usetikzlibrary{decorations.pathreplacing}

\newtheorem{theorem}{Theorem}

\newcommand{\btilde}{\widetilde{b}}
\newcommand{\qtilde}{\widetilde{q}}
\newcommand{\ftilde}{\widetilde{f}}

\DeclareMathOperator{\dof}{dof}

\title{Turing patterns in a 3D morpho-chemical bulk-surface reaction-diffusion system for battery modeling}

\author{Massimo Frittelli\thanks{University of Salento, Department of Mathematics and Physics ``E. De Giorgi'', Via per Arnesano, 73100 Lecce, Italy. Emails: \texttt{massimo.frittelli@unisalento.it}, \texttt{ivonne.sgura@unisalento.it}}, Ivonne Sgura\samethanks, Benedetto Bozzini\thanks{Department of Energy, Politecnico di Milano, Milano, Italy. Email: \texttt{benedetto.bozzini@polimi.it}}}

\date{}

\newcommand*\samethanks[1][\value{footnote}]{\footnotemark[#1]}

\begin{document}

\maketitle

\begin{abstract}
In this paper we introduce a bulk-surface reaction-diffusion (BSRD) model in three space dimensions that extends the DIB morphochemical model to account for the electrolyte contribution in the application, in order  to study structure formation during discharge-charge processes in batteries. Here we propose to approximate the model by the Bulk-Surface Virtual Element Method on a tailor-made mesh that proves to be competitive with fast bespoke methods for PDEs on Cartesian grids. We present a selection of numerical simulations that accurately match the classical morphologies found in experiments. Finally, we compare the Turing patterns obtained by the coupled 3D BS-DIB model with those obtained with the original 2D version.
\end{abstract}

\section*{Keywords}
Batteries; Metal electrode; Electrodeposition; Bulk-Surface Reaction-Diffusion Systems; Bulk-Surface Virtual Element Method; Turing Patterns

\section*{MSC 2020}
65M60; 65M50; 65N40; 65P40

\section{Introduction}
The formation of spatio-temporal structures in electrodeposition is a relevant physical phenomenon, as it impacts several applications, ranging from the durability and efficiency of batteries to electroplating \cite{bozzini2013spatio}.  The onset of spatio-temporal structures on the cathodic surface was proven to be initiated by a Turing morphogenetic mechanism,  where the physics are modeled by a suitable reaction-diffusion system (RDS), known as DIB model,  whose spatial domain is the electrodic surface \cite{bozzini2013spatio}.  In the DIB model, the spatial domain is assumed to be fixed and does not change over time, as the growth/corrosion effects are fully modeled by the dynamics of the system.  By tweaking the parameters of the DIB model, it is possible to successfully simulate spatial \cite{bozzini2013spatio} or spatio-temporal patterns \cite{lacitignola2015} of various morphological classes that are experimentally observed under appropriate physical and chemical conditions; these include spatial patterns such as spots,  holes,  stripes, labyrinths, and spiral waves. The effectiveness of the DIB model has justified the development of extensions and ameliorations, such as the introduction of cross-diffusion \cite{DIB_cross_diff} and the generalization of the spatial domain to be a curved surface \cite{DIB_sphere}.

As it stands, one of the limitations of the DIB model is that it does not fully accounts for the effects of non-uniform electrolite concentration in a neighborhood of the electrode.  Experimentally, such non uniform concentration is induced by the spatial structures arising on the electrode and, in turn, affects further structure development.  In this regard, the electrode-electrolyte system has a two-way coupling that, in the long run, can drastically affect the resulting morphological class. In this paper we propose the \emph{bulk-surface DIB (BS-DIB) model} in three space dimensions to fill this gap.  In the proposed model, the surface represents the electrode (where the electrodeposition takes place), while the 3D bulk models the electrolite. The physical two-way coupling mentioned above causes the proposed model to take the form of a \emph{coupled bulk-surface reaction-diffusion system (BS-RDS)} \cite{madzvamuse_bsrds}.

For domains of general shape, different numerical methods were developed for the spatial approximation of BS-RDSs,  such as the Bulk-Surface Finite Element Method (BS-FEM) \cite{bsfem}, the Cut Finite Element Method \cite{cutFEM}, unfitted finite element methods \cite{unfitted_bsfem}, and meshless kernel methods \cite{kernel_bspde}, just to mention a few. In all of these methods, the spatially discrete problem takes the form of a large ODE system, whose dimension is equal to the number of spatial degrees of freedom.  Thus, the high level of spatial resolution required by RDSs and BS-RDSs, together with the curse of dimensionality (3D) makes the numerical approximation of the BS-DIB a challenging computational task.  In the present context, where the bulk domain is a cube, a bespoke tensorized technique called Matrix-Oriented Finite Element Method (MO-FEM) \cite{mofem3d} can be exploited to take advantage of the special geometry and drastically reduce the computational times. However, it is worth noting that the BS-DIB model exhibits spatial patterns only in a neighborhood of the surface, hence a uniform spatial discretization is computationally inefficient. For this reason, we exploit the geometric flexibility of the Bulk-Surface Virtual Element Method (BS-VEM) \cite{bsvem_parabolic} to adopt a graded cubic mesh that is highly refined close to the surface and much coarser away from the surface.  Such a mesh is simultaneously graded and entirely composed of cubic-shaped elements. Such a combination entails the presence of hanging nodes and edges, that are naturally handled by the BS-VEM and are not admissible in the BS-FEM. Compared to the MO-FEM, the BS-VEM on such mesh exhibits shorter computational times on equal level of spatial refinement in a neighborhood of the surface -where high spatial accuracy is actually required- and produce patterns of the same morphological class.  It needs to be noted that Turing patterns are highly sensitive to initial conditions, which are bound to be different between MO-FEM ad BS-VEM since the spatial meshes are different, hence obtaining the same morphological class with both  methods is a sensible benchmark.

We present a wide range of numerical simulations for both the 2D BIB and the 3D BS-DIB models in equal parameters, to showcase the effect of bulk-surface coupling.  From the experiments we draw the following conclusions.  First, the BS-DIB model appears to have a large Turing region compared to the DIB model. In fact, for several choices of the parameters outside the Turing region of the DIB model, only the BS-DIB model exhibits spatial patterns. Second, when the DIB model exhibits spatial patterns, the BS-DIB still exhibits patterns, but if different morphological class,  thereby further highlighting the impact of the bulk-surface coupling. A rigorous analysis of the Turing instability for the BS-DIB model will be addressed in future work.

The structure of the paper is as follows. In Section \ref{sec:da_model} we introduce the BS-DIB model, we give its physical interpretation and we analyse the stability of a relevant spatially uniform equilibrium in the absence of diffusion. In Section \ref{sec:bsvem} we recall the BS-VEM and we present a bespoke graded polyhedral mesh that allows the BS-VEM to outperform the MO-FEM when solving the BS-DIB model.  In Section \ref{sec:simulations_turing} we present an extensive list of numerical experiments that empirically demonstrate the effect of the bulk-surface coupling on pattern formation. Conclusions are drawn in Section \ref{sec:conclusions}.

\section{The bulk-surface DIB model on the cube}
\label{sec:da_model}
Let $\Omega = [0,L]^3$ be a cube of edge $L>0$, let $\Gamma := [0,L]^2\times \{0\}$ be the bottom face of $\Gamma$, let $\Gamma_T = [0,L]^2\times \{L\}$ be the top face of $\Omega$ and let $\Gamma_L = \partial \Omega \setminus (\Gamma\cup\Gamma_T)$ be the union of the lateral faces of $\Omega$.  Let $T>0$ be the final time.  
The bulk-surface DIB model seeks to find four functions $b,q:\Omega\times [0,T] \rightarrow\mathbb{R}$ and $\eta,\theta:\Gamma\times [0,T] \rightarrow\mathbb{R}$ that fulfil the following system of four PDEs:
\begin{equation}
\label{model}
\begin{cases}
\dot{b} - \Delta b = f_1(b) \qquad \text{in}\ \Omega;\\
\dot{q} - d_\Omega \Delta q = f_2(q) \qquad \text{in}\ \Omega;\\
\dot{\eta} - \Delta_\Gamma \eta = f_3(b,\eta,\theta) \qquad \text{on}\ \Gamma;\\
\dot{\theta} - d_\Gamma\Delta_\Gamma \theta = f_4(q,\eta,\theta) \qquad \text{on}\ \Gamma,
\end{cases}
\end{equation}
complemented with the following boundary conditions for the bulk variables $b$ and $q$
\begin{equation}
\label{boundary_conditions}
\begin{cases}
\nabla b \cdot \mathbf n = - f_3(b,\eta,\theta)\psi_\eta \qquad \text{on} \ \Gamma;\\
\nabla q \cdot \mathbf n = - f_4(q,\eta,\theta)\psi_\theta\qquad \text{on} \ \Gamma;\\
\nabla b \cdot \mathbf n = 0 \qquad \text{on} \ \Gamma_L;\\
\nabla q \cdot \mathbf n = 0 \qquad \text{on} \ \Gamma_L;\\
b = b_0 \qquad \text{on}\ \Gamma_T;\\
q = q_0 \qquad \text{on}\ \Gamma_T,
\end{cases}
\end{equation}
and the following boundary conditions for the surface variables $\eta$ and $\theta$:
\begin{equation}
\label{boundary_conditions_surf}
\begin{cases}
\nabla \eta \cdot \mathbf n = 0 \qquad \text{on}\ \partial \Gamma;\\
\nabla \theta \cdot \mathbf n = 0 \qquad \text{on}\ \partial \Gamma.
\end{cases}
\end{equation}
In \eqref{model}, $\Delta$ is the Laplace operator in $\Omega$, while $\Delta_\Gamma$ is the Laplace-Beltrami operator on $\Gamma$ (which coincides with the two-dimensional Laplacian since $\Gamma$ is flat), $d_\Omega$ and $d_\Gamma$ are diffusion coefficients, and $b_0,q_0 \in\mathbb{R}$. 
In the equations for the bulk species $b$ and $q$, the kinetics $f_1,f_2$ are defined as follows:
\begin{align}
\label{kinetics_b}
&f_1(b) := -k_b(b-b_0);\\
\label{kinetics_q}
&f_2(q) := -k_q(q-q_0);
\end{align}
$b$ represents the concentration of the electroactive cations (precursors of metal that is electrodeposited during the recharge cycle), present exclusively in the bulk.
$q$ represents the bulk concentration of an additive species that is adsorbed at the cathode as a way of controlling shape change with the coverage degree expressed by the variable $\theta$. 
$b_0=b_{bulk}$ and $q_0=q_{bulk}$ represent the “bulk concentrations”, that prevail at equilibrium, when the bulk is homogeneous. The physical meaning of the terms $(b-b_0)$ and $(q-q_0)$, with $k_b,k_q>0$ is first-order homogeneous reaction kinetics describing the tendency of the reagent to re-establish the equilibrium concentration. This can be considered a very simple model of a situation in which $b$, $q$  are the concentrations of the species involved in the electrodic reaction in the electroactive form, that is generated by the decomposition of some precursor (e.g.metallic ion with a ligand that keeps the ion in solution in non-electroactive form, from which the electroactive species forms by decomposition of the complexed one): a lower-than-equilibrium local concentration of $b$, $q$ (e.g. by cathodic consumption) generates new $b, q$ (by decomposition of the complexed form); while a higher-than-equilibrium concentration (e.g. by anodic injection) generates a consumption (e.g. by reaction with the ligand, yielding the non-electroactive form).
In the surface equations for the surface species $\eta$ and $\theta$, the kinetics $f_3$ and $f_4$ are defined as follows:
\begin{align}
\label{kinetics_eta}
&f_3(b,\eta,\theta) :=  \rho [A_1b(1-\theta)\eta - A_2\eta^3 - B(\theta-\alpha)];\\
\label{kinetics_theta}
&f_4(q,\eta,\theta) := \rho \left[q(1+k_2\eta)(1-\theta)(1-\gamma(1-\theta)) -\frac{D}{C}(1+k_3\eta)\theta (1+\gamma\theta)\right],
\end{align}
and correspond to the modified DIB model source terms accounting for the bulk contributions, as explained below, to account for the electrolyte physics.

The physical meaning of $f_3$ is that the Butler-Volmer type electrokinetic term   is scaled by the concentration of the electroactive species at the surface $b_|\Gamma$. This corresponds to first-order phenomenological kinetics and it includes naturally mass-trasport effects from the bulk to the surface (i.e. the mass-transport of the electroactive species $b$ present in the bulk and reacting electrochemically at the surface). It is worth noting that the term $-A_2\eta^3$ in the 2D form of DIB represents collectively all hindrances to $\eta$ (growth) resulting from the establishment of high values of $\eta$.  The most straightforward interpretation of such hindrances is mass-transport limitation that – in a Gileadi-type framework – can be approximately accounted for with a negative cubic correction to the I-V curve. As expounded above, in the bulk-surface context mass-transport limitations can be naturally accounted for by multiplying the term linear in $\eta$ by the surface concentration $b_|\Gamma$ of the reactant $b$ present in the bulk and describing the electroactive species. Nevertheless, in the bulk-surface version of DIB it is worth retaining the cubic term in $\eta$, because this can account for other hindrances to metal growth (i.e. beyond mass-transport from the bulk to the reactive surface) appearing at high metal plating rates, such as cathodic passivation. 

The physical meaning of $f_4$ is that - coherently with the Langmuir adsorption model with monomolecular adsorption reaction - the adsorption term is directly proportional to the amount of the bulk species. In the case of a heterogeneous bulk phase, the relevant value of the bulk form is that in a neighborhood of the surface (commonly, referred to as the ``catholyte''): $q_|\Gamma$. In \eqref{kinetics_b}-\eqref{kinetics_q}-\eqref{kinetics_eta}-\eqref{kinetics_theta},  $d_\Omega, d_\Gamma,\psi_\eta,\psi_\theta, b_0,q_0,k_b,k_q,k_2,k_3, \rho,\alpha,\gamma, A_1,A_2,B,C,D$ are positive coupling parameters. 

The coupling BCs for the bulk equations at the interface between the ``growing surface $\Gamma$'' and the ``bulk Q'' can thus naturally be written as in $(2)_,  \ (2)_2$ indicating that the flux of bulk species to the surface is opposite to their consumption rates at the surface. More specifically, the physical meaning of the form of Eq. $(2)_, \ (2)_2$ is that, even though $b$ and $q$ coming from the bulk are consumed at the interface (i.e. in correspondence of their values $b_|\Gamma$  and $q_|\Gamma$ ) to yield $\eta$ and $\theta$ only in a specific term of $f_3$ and $f_4$,  the negative terms of these equations have the effect of injecting $b$ and $q$ into the bulk. Cases in which $\theta > \alpha$ (i.e. the adsorbate enhances electrodeposition, e.g. by resonant tunnelling effects), can be regarded as a special case of interfacial $b$ consumption, accounted for through the BCs. Thus the net formation rate of $\eta$ and $\theta$ is proportional to the fluxes of $b$ and $q$ to the surface.

 $\nabla b \cdot \mathbf n, \nabla q \cdot \mathbf n$ denote the gradients normal to the boundary (surface) $\Gamma$, while $\psi_\eta, \psi_\theta$  are constants the role of which is to adjust the dimensionality of the equations and the physical meaning of which is explained below. $\psi_\eta$  converts adatoms (i.e. the surface species generated by the reaction (consumption) of $b$ at the surface) into morphological units (the quantities actually described by $\eta$). Thus one can write: $\eta=\psi_\eta b$  and $\psi_\eta$  can be regarded as a constant as far as the density of morphological units (ca. step) is proportional to the adatom density.   $\psi_\theta$ expresses an isotherm, since it connects a bulk concentration ($q$) into a surface density ($\theta$), hence, in as far as the isotherm can be linearised, we can write: $\theta = \psi_\theta q$.

The BCs for $b$ and $q$ are:\\
i)  non-linear coupling BCs on the surface $\Gamma$, implying coupling with $\eta$ and $\theta$; ii) Dirichlet BCs on the face of $\Omega$ opposite to $\Gamma$ that is located far enough from the “bottom face” where reaction takes place, so that concentration gradients induced by reactivity have died out here; iii) zero Neumann BCs  (zero flux) on the residual faces of $\Omega$ (see Figure 1). The values of the bulk variables are thus set to their equilibrium values $b_0,q_0$ respectively.
The BCs for the modifed DIB model on $\partial \Gamma$ are also zero Neumann BCs.
 
\begin{figure}[h!]
\begin{center}
{\includegraphics[scale=.45]{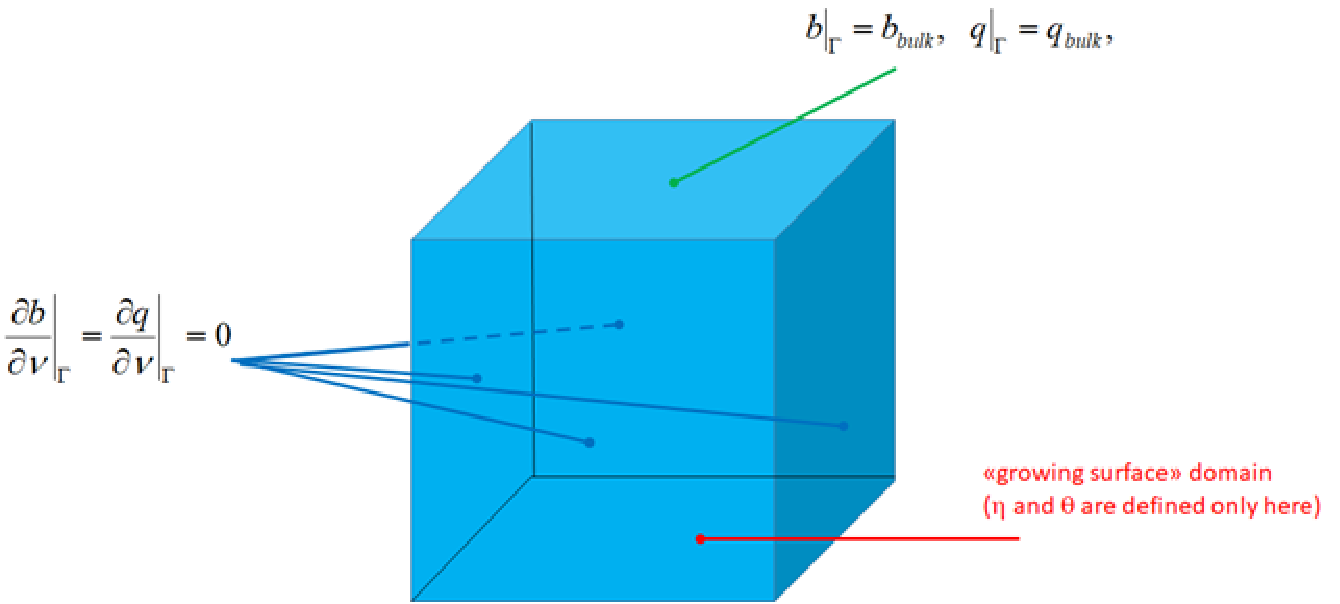}}
\caption{Domain geometry and BCs for the bulk variables $b$ and $q$.}
\end{center}
\label{fig:BSdomain}
\end{figure}

\noindent
If $b_0=q_0 = 1$, the following is a spatially homogeneous equilibrium for the system \eqref{model}-\eqref{boundary_conditions}:
\begin{equation}
\label{equilibrium}
(b^*,q^*, \eta^*, \theta^*) = (1,1,0,\alpha).
\end{equation}
The initial conditions are prescribed as follows:
\begin{align}
&b(\mathbf x,0) = b_0;\\
&q(\mathbf x, 0) = q_0;\\
&\eta(\mathbf x, 0) = r_\eta(\mathbf x);\\
&\theta(\mathbf x, 0) = r_\theta(\mathbf x),
\end{align}
where $r_\eta$ and $r_\theta$ are random spatial data that fulfil
\begin{align}
&r_\eta(\mathbf x) \in \eta_e+[0, 1e-2]\qquad \forall \mathbf x \in \Gamma;\\
&r_\theta(\mathbf x) \in \theta_e+[-1e-2, 1e-2] \qquad \forall \mathbf x \in \Gamma.
\end{align}
The model \eqref{model}-\eqref{boundary_conditions} is a generalisation of the DIB model in 2D,  which takes the form
\begin{equation}
\label{model2d}
\begin{cases}
\dot{\eta} - \Delta_\Gamma \eta = f(\eta,\theta) \qquad \text{in}\ Q;\\
\dot{\theta} - d_\theta\Delta_\Gamma \theta = g(\eta,\theta) \qquad \text{in}\ Q;\\
\nabla \eta \cdot \mathbf n = 0 \qquad \text{on}\ \partial Q;\\
\nabla \theta \cdot \mathbf n = 0 \qquad \text{on}\ \partial Q,
\end{cases}
\end{equation}
where $Q\subset\mathbb{R}^2$ is a compact 2D domain,  while $f(\eta,\theta) = f_3(1,\eta,\theta)$ and $g(\eta,\theta) = f_4^\gamma(1,\eta,\theta)$, see \cite{bozzini2013spatio}.  After \cite{lacitignola2015}, the parameter values are fixed as follows:
\begin{align}
\label{parameter_values}
\begin{cases}
& d_\Omega = 1;\\
& d_\Gamma = 20;\\
&k_2 = 2.5;\\
&k_3 = 1.5;\\
&\rho = 1;\\
&\alpha = 0.5;\\
&A_1 = 10;\\
&D = \frac{C(1-\alpha)(1-\gamma+\gamma\alpha)}{(\alpha(1+\gamma\alpha))},
\end{cases}
\end{align}
while $(B,C)$ are bifurcation parameters.
The novel parameters instead (the ones that do not appear in the 2D DIB model \eqref{model2d}) are fixed as follows:
\begin{align}
\label{parameters_all_experiments}
\begin{cases}
&b_0 = 1;\\
&q_0 = 1;\\
&k_b = 1;\\
&k_q = 1.
\end{cases}
\end{align}
Finally, the other parameters will be changed in a representative series of cases that will be discussed below. 
If $F=(f_1,f_2,f_3,f_4)$ and $\xi=(b, q,\eta,\theta)^T$, it is easy to show that 
\begin{equation}
\label{eq}
F(\xi*)=0     \quad \text{if}   \quad  \xi*=(b^*, q^*, \eta^*,\theta^*) = (b_0, q_0, 0, \alpha),
\end{equation}
when in the BS-DIB model all parameters different from (C, B) are kept fixed at the previous values. Note that the part of this equilibrium concerning the variables $\eta,\theta$ coincides with the equilibrium of the original DIB (2D) model. 
In the next section we shall analyze the stability of such equilibrium in the special case $\gamma =0$,  when the condition on the parameter $D$ in \eqref{parameter_values} boils down to
\begin{equation}
\label{gamma_0}
\gamma=0, \ D=q_0 \frac{C(1-\alpha)}{\alpha}.
\end{equation}

\subsection{Stability in absence of diffusion}
To study the arising of the diffusion-driven or Turing instability in the BS-DIB model, it is necessary to prove that the equilibrium \eqref{eq} is stable in absence of diffusion. In fact, the model \eqref{model},  deprived of diffusion and linearized around the equilibrium \eqref{equilibrium} is
\begin{equation}
\boldsymbol \xi_t = J(\xi^*) \boldsymbol (\xi-\xi^*),
\end{equation}
where $\boldsymbol \xi = (b, q, \eta, \theta)^T$ and $J$ is the Jacobian of the kinetics evaluated at the equilibrium \eqref{eq}. The matrix $J$ has the following structure
\begin{equation}
J = \begin{pmatrix}
J_{\Omega} & 0\\
J_h & J_\Gamma
\end{pmatrix},
\end{equation}
where, if $\gamma = 0$, the blocks of $J$ are as follows:
\begin{align}
&J_\Omega :=
\begin{pmatrix}
f_{1,b} & f_{1,q}\\
f_{2,b} & f_{2,q}
\end{pmatrix} = 
\begin{pmatrix}
- k_b & 0\\
0 & -k_q
\end{pmatrix};\\
&J_h := \rho
\begin{pmatrix}
f_{3,b} & f_{3,q}\\
f_{4,b} & f_{4,q}
\end{pmatrix} = 
\rho\begin{pmatrix}
A_1 (1-\theta)\eta & 0\\
0 & C(1+k_2\eta)(1-\theta)
\end{pmatrix};\\
&J_\Gamma := 
\begin{pmatrix}
f_{3,\eta} & f_{3,\theta}\\
f_{4,\eta} & f_{4,\theta}
\end{pmatrix} = 
\rho\begin{pmatrix}
b_0 A_1 (1-\alpha) & -B\\
q_0 C(k_2-k_3)(1-\alpha) & -\frac{q_0 C}{\alpha}
\end{pmatrix}.
\end{align}
Thanks to the diagonal structure of $J_\Omega$ it holds that
\begin{equation}
\label{eigenvalues_absence_diffusion}
\det(J - \lambda I) = (\lambda + k_b)(\lambda + k_q)\det(J_\Gamma - \lambda I).
\end{equation}
It follows that two eigenvalues of $J$ are $\lambda_1 = -k_b < 0$ and $\lambda_2 = -k_q < 0$.  We are left to determine when the eigenvalues of $J_\Gamma$ are negative.  This happens if and only if $Trace J_\Gamma < 0$ and $\det J_\Gamma > 0$. Now:
\begin{equation}
Trace J_\Gamma = \rho\left(b_0 A_1(1-\alpha) - \frac{q_0 C}{\alpha}\right) < 0 \Longleftrightarrow C > \frac{b_0}{q_0}A_1\alpha(1-\alpha),
\end{equation}
and
\begin{equation}
\det J_\Gamma = \rho^2 \Big(BC q_0 (k_2-k_3)(1-\alpha) - A_1C  b_0 q_0 \frac{1-\alpha}{\alpha}\Big) > 0 \Longleftrightarrow B > \frac{A_1 b_0}{\alpha(k_2-k_3)}.
\end{equation}
We obtain the following result.
\begin{theorem}
If $\gamma = 0$, the equilibrium \eqref{equilibrium} is stable in the absence of diffusion if and only if
\begin{equation}
\label{stability_equilibrium_no_diffusion}
B > \frac{A_1 b_0}{\alpha(k_2-k_3)} \qquad \wedge \qquad C > \frac{b_0}{q_0}A_1\alpha(1-\alpha).
\end{equation}
In addition, if $A_1$,  $k_2$, $k_3$, $\alpha$ are as in \eqref{parameter_values},  the condition \eqref{stability_equilibrium_no_diffusion} specializes to
\begin{equation}
\label{stability_equilibrium_no_diffusion_specialized}
B > 20 b_0 \qquad \wedge \qquad C > 2.5 \frac{b_0}{q_0}.
\end{equation}
Furthermore, if $b_0 = q_0 = 1$,  the condition  \eqref{stability_equilibrium_no_diffusion_specialized} further reduces to
\begin{equation}
B > 20 \qquad \wedge \qquad C > 2.5.
\end{equation}
\end{theorem}

Of course, as far as the conditions for Turing pattern formation are concerned, a specific study has to be carried out considering the full jacobian J of Eq. (10). This is an important study topic in its own right, that - nevertheless - has no impact on the analysis presented in this reserach.  Since treating this problem exhaustively would be beyond the scope of the present work,  we leave to a subsequent publication.
Our present approach is thus to solve numerically the BS-DIB model for a representative selection of parameter couples $(B,C)$ generating the whole set of Turing pattern morphologies, as described in \cite{amospaper, mgq}. Moreover, we shall will fix the model parameters in the bulk as in \eqref{parameters_all_experiments}.  In this scenario, our aim is to tune the coupling parameters $\psi_\eta, \psi_\theta$ to study from the numerical point of view the effect of coupling with the the bulk on the morphological structure of the Turing patterns in the classes studied in \cite{amospaper}.
We recall that the numerical approximation reaction-diffusion systems in 3D is not straightforward because the pattern requires a very fine 3D mesh that provides sufficient spatial resolution and a longtime integration to reach the asymptotic steady state.   For this reason, we shall apply the BS-VEM method studied in \cite{bsvem_parabolic}, for the space discretization of the 3D domain and surface and the IMEX Euler method as time solver.  Hence, in the 3D case,  if the cubic domain is approximated with a Cartesian grid,  at least a million of unknowns at each time iteration are required.  The usual implementation will thus end up to a sequence of linear systems where the coefficient matrix for each species is sparse, but prohibitively large. An efficient strategy to deal with this issue is the matrix-oriented (MO) approach \cite{mofem3d} where, thanks to the Cartesian structure of the numerical grid,  the fully discrete problem is transformed to a sequence of Sylvester matrix equations, that are solved in the spectral space.  However,  since BS-DIB model produces spatial patterns only in a neighborhood of the surface $\Gamma$, we devise a tailor-made graded polyhedral bulk-surface mesh where the BS-VEM proves to be a competitive alternative, since such a graded mesh avoids unnecessary refinement (and degrees of freedom) away from the surface $\Gamma$.  One of the major advantages of this fact is that the BS-VEM, thanks to the flexibility of polyhedral meshes, can still be used used on domains of general shape, where MO techniques might not apply.

\section{The Bulk Surface Virtual Element Method for the BS-DIB model}
\label{sec:bsvem}

Formulating a Bulk-Surface Virtual Element Method (BSVEM) \cite{bsvem_parabolic} for the model \eqref{model} requires several steps. We start by rewriting the model \eqref{model} in such a way that the boundary conditions lend themselves to a BSVEM discretization.

\subsection{Step 1: Rewriting the model with homogeneous boundary conditions}
In the presence of non-zero Dirichlet boundary conditions, it is well-known \cite{quarteroni_book} that it is first necessary to rewrite the PDE problem in such a way that the Dirichlet conditions are homogeneous. To this end,  we define the following auxiliary variables and kinetics:
\begin{align}
&\btilde := b - b_0;\\
&\qtilde := q - q_0;\\
&\ftilde_1(\btilde) := f_1(\btilde + b_0);\\
&\ftilde_2(\qtilde) := f_2(\qtilde + q_0);\\
&\ftilde_3(\btilde, \eta,\theta) := f_3(\btilde + b_0, \eta, \theta);\\
&\ftilde_4(\qtilde, \eta,\theta) := f_4(\qtilde + q_0, \eta, \theta).
\end{align}
With the above definitions, the model becomes
\begin{equation}
\label{model_auxiliary}
\begin{cases}
\dot{\btilde} - \Delta \btilde = \ftilde_1(\btilde) \qquad \text{in}\ \Omega;\\
\dot{\qtilde} - d_\Omega \Delta \qtilde = \ftilde_2(\qtilde) \qquad \text{in}\ \Omega;\\
\dot{\eta} - \Delta_\Gamma \eta = \ftilde_3(\btilde,\eta,\theta) \qquad \text{on}\ \Gamma;\\
\dot{\theta} - d_\Gamma\Delta_\Gamma \theta = \ftilde_4(\qtilde,\eta,\theta) \qquad \text{on}\ \Gamma,
\end{cases}
\end{equation}
which, this time, is conveniently endowed with completely homogeneous boundary conditions:
\begin{equation}
\label{boundary_conditions_auxiliary}
\begin{cases}
\nabla \btilde \cdot \mathbf n = - \ftilde_3(\btilde,\eta,\theta)\psi_\eta \qquad \text{on} \ \Gamma;\\
\nabla \qtilde \cdot \mathbf n = - \ftilde_4(\qtilde,\eta,\theta)\psi_\theta\qquad \text{on} \ \Gamma;\\
\nabla \btilde \cdot \mathbf n = 0 \qquad \text{on} \ \Gamma_L;\\
\nabla \qtilde \cdot \mathbf n = 0 \qquad \text{on} \ \Gamma_L;\\
\btilde = 0 \qquad \text{on}\ \Gamma_T;\\
\qtilde = 0 \qquad \text{on}\ \Gamma_T,
\end{cases}
\end{equation}

\subsection{Step 2: Weak formulation}
To write a discrete formulation of the auxiliary problem \eqref{model_auxiliary}-\eqref{boundary_conditions_auxiliary}, we define the space of trivariate spatial functions that ensure the well-posedness of the model \eqref{model_auxiliary} and fulfil the boundary conditions \eqref{boundary_conditions_auxiliary}:
\begin{equation}
H^1_B(\Omega) := \{u \in H^1(\Omega) \ | \ u_{|\Gamma_T} = 0 \ \wedge\  u_{|\Gamma}\ \in  H^1(\Gamma)\}.
\end{equation}
The dual space of $H^1_B(\Omega)$ will be denoted by $H^{-1}_B(\Omega)$.  Following \cite{bsvem_parabolic}, the weak formulation of \eqref{model_auxiliary}-\eqref{boundary_conditions_auxiliary} is: find $\btilde,\qtilde \in L^2([0,T]; H^1_B(\Omega))$ and $\eta,\theta \in L^2([0,T]; H^1(\Gamma))$ with $\dot{\btilde}, \dot{\qtilde} \in L^2([0,T]; H^{-1}_B(\Omega))$ and $\dot{\eta}, \dot{\theta} \in L^2([0,T]; H^{-1}(\Gamma))$ such that
\begin{equation}
\label{weak_formulation}
\begin{cases}
\displaystyle\int_\Omega \dot{\btilde}\varphi + \int_\Omega \nabla \btilde \cdot \nabla \varphi = \int_\Omega \ftilde_1(\btilde)\varphi - \psi_\eta\int_\Gamma  \ftilde_3(\btilde, \eta, \theta)\varphi;\\
\displaystyle\int_\Omega \dot{\qtilde}\varphi + d_\Omega\int_\Omega \nabla \qtilde \cdot \nabla \varphi = \int_\Omega \ftilde_2(\qtilde)\varphi - d_\Omega\psi_\theta \int_\Gamma  \ftilde_4(\qtilde, \eta, \theta)\varphi;\\
\displaystyle\int_\Gamma \dot{\eta}\psi + \int_\Gamma \nabla_\Gamma \eta \cdot \nabla_\Gamma \psi = \int_\Gamma \ftilde_3(\btilde,\eta, \theta)\psi;\\
\displaystyle\int_\Gamma \dot{\theta}\psi + d_\Gamma \int_\Gamma \nabla_\Gamma \theta \cdot \nabla_\Gamma \psi = \int_\Gamma \ftilde_4(\qtilde,\eta, \theta)\psi,
\end{cases}
\end{equation}
for all $\varphi \in L^2([0,T]; H^1_B(\Omega))$ and $\psi \in L^2([0,T]; H^1(\Gamma))$.

\subsection{Step 3: Spatially discrete formulation}
We will now describe the spatial formulation obtained by the BSVEM following \cite{bsvem_parabolic}. The choice of the BSVEM to solve the model \eqref{model} is motivated by the possibility of using graded meshes that make the BSVEM particularly competitive in this case by avoiding unnecessarily refinement away from the surface $\Gamma$.  The choice of a convenient mesh will be illustrated in the next Section. For now, we illustrate the BSVEM for arbitrary meshes. 

Let us decompose the bulk $\Omega$ as the union of non-overlapping polyhedra, $\Omega = \cup_{E\in\mathcal{E}_h} E$.  If $\mathcal{F}_f$ is the set of the faces of $\mathcal{E}_h$ that are contained in $\Gamma$, then we can write $\Gamma = \cup_{F\in\mathcal{F}_h} F$.  For a face $F\in\mathcal{F}_h$, the \emph{boundary space} $\mathbb{B}(\partial F)$ is defined by
\begin{equation}
\mathbb{B}(\partial F) := \{v \in\mathcal{C}^0(\partial F) \ | \ v_e \in\mathbb{P}_1(e) \ \forall e \in\text{edges}(F)\}.
\end{equation}
The \emph{preliminary space} of a face $F$ is defined by
\begin{equation}
\widetilde{\mathbb{V}}(F) := \{v \in H^1(F) \ | \ v_{|\partial F} \in \mathbb{B}(\partial F) \wedge \Delta v \in\mathbb{P}_1(F)\}.
\end{equation}
The $H^1$ projector on faces $\Pi^\nabla_F: \widetilde{\mathbb{V}}(F) \rightarrow \mathbb{P}_1(F)$ is defined, for any $v \in \widetilde{\mathbb{V}}(F)$ by
\begin{equation}
\int_F \nabla (v - \Pi^\nabla_F v) \cdot \nabla p = 0 \ \forall p\in\mathbb{P}_1(F) \wedge \int_F (v-\Pi^\nabla_F v) = 0. 
\end{equation}
Then, the \emph{enhanced VEM space} on the face $F$ is defined by
\begin{equation}
\mathbb{V}(F) := \left\{v \in \widetilde{\mathbb{V}}(F)\  \middle| \ \int_F (v - \Pi^\nabla_F v)p = 0 \ \forall p\in\mathbb{P}_1(F)\right\}.
\end{equation}
For a polyhedron $E\in\mathcal{E}_h$,  the \emph{boundary space} $\mathbb{B}(\partial E)$ is defined by
\begin{equation}
\mathbb{B}(\partial E) := \{u \in  \mathcal{C}^0(\partial E) \ |\  u_{|F}\in\mathbb{V}(F)\ \forall F \in \text{faces}(E)\}.
\end{equation}
At this point, the \emph{preliminary VEM space} on $E$ is defined by
\begin{equation}
\widetilde{\mathbb{V}}(E) := \{u \in H^1(E) \ | \ u_{|\partial E} \in \mathbb{B}(\partial E) \wedge \Delta u \in \mathbb{P}_1(E)\}.
\end{equation} 
The $H^1$ projector $\Pi^\nabla_E: \widetilde{\mathbb{V}}(E) \rightarrow \mathbb{P}_1(E)$ on the polyhedron $E$ is defined, for each $u\in \widetilde{\mathbb{V}}(E)$, by
\begin{equation}
\int_E \nabla (u - \Pi^\nabla_E u) \cdot \nabla p = 0 \ \forall p\in\mathbb{P}_1(E) \wedge \int_E (u-\Pi^\nabla_E u) = 0. 
\end{equation} 
Finally, the \emph{enhanced VEM space} on the polyhedron $E$ is defined by
\begin{equation}
\mathbb{V}(E) := \left\{u \in \widetilde{\mathbb{V}}(E)\  \middle| \ \int_E (u - \Pi^\nabla_E u)p = 0 \ \forall p\in\mathbb{P}_1(E)\right\}.
\end{equation}
It is well-known that the degrees of freedom in $\mathbb{V}(F)$ and $\mathbb{V}(E)$ are the pointwise values on vertexes, see \cite{ahmad}. The global VEM spaces are defined by matching the degrees of freedom across elements. To this end,  let $\mathbb{S}_\Gamma$ and $\mathbb{S}_\Omega$ be the 1-skeleton of $\Gamma$ and the 2-skeleton of $\Omega$, respectively, defined by
\begin{align}
\mathbb{S}_\Gamma := \bigcup_{F \in\mathcal{F}_h} \partial F, \qquad \mathbb{S}_\Omega := \bigcup_{E \in \mathcal{E}_h} \partial E.
\end{align}
The global VEM spaces $\mathbb{V}_\Gamma$ and $\mathbb{V}_\Omega$ are then defined as
\begin{align}
&\mathbb{V}_\Gamma := \{v\in H^1(\Gamma) \ |\ v \in\mathcal{C}^0(\mathbb{S}_\Gamma) \ \wedge \ v_{|F} \in \mathbb{V}(F) \ \forall F \in \mathcal{F}_h\};\\
&\mathbb{V}_\Omega := \{u \in H^1(\Omega) \ | \ u\in\mathcal{C}^0(\mathbb{S}_\Omega) \ \wedge \ u_{|E}\in\mathbb{V}(E) \ \forall E \in \mathcal{E}_h \ \wedge \ u(x,y,L) = 0 \ \forall (x,y) \in [0,L]^2\}.
\end{align}
Notice that the space $\mathbb{V}_\Omega$ reflects the homogeneous Dirichlet boundary conditions of the continuous counterpart $H^1_B(\Omega)$. To obtain a spatially discrete counterpart of the weak formulation \eqref{weak_formulation}, we need suitable discrete bilinear forms. Following \cite{bsvem_parabolic},  for all $F\in\mathcal{F}_h$, $E\in\mathcal{E}_h$, $v,w \in \mathbb{V}(F)$ and $u,z\in\mathbb{V}(E)$, we define
\begin{align}
&m_F(v,w) := \int_F \Pi^0_F v \Pi^0_F w + h_F^2 \langle \dof(v-\Pi^0_F v), \dof(w-\Pi^0_F w) \rangle;\\
&a_F(v,w) := \int_F \nabla\Pi^\nabla_F v \cdot \nabla \Pi^\nabla_F w + \langle \dof(v-\Pi^0_F v), \dof(w-\Pi^0_F w) \rangle;\\
&m_E(u,z) := \int_E \Pi^0_E u \Pi^0_E z + h_E^3 \langle \dof(u-\Pi^0_E u), \dof(z-\Pi^0_E z) \rangle;\\
&a_E(u,z) := \int_E \nabla\Pi^\nabla_E u \cdot \nabla \Pi^\nabla_E z + h_E\langle \dof(u-\Pi^0_E u), \dof(z-\Pi^0_E z) \rangle,
\end{align}
where $h_F$ and $h_E$ are the diameters of $F$ and $E$, respectively. Let $m_h^\Gamma, a_h^\Gamma : \mathbb{V}_\Gamma \times \mathbb{V}_\Gamma \rightarrow\mathbb{R}$ and $m_h^\Omega, a_h^\Omega : \mathbb{V}_\Omega \times \mathbb{V}_\Omega \rightarrow\mathbb{R}$ be the corresponding global forms. Furthermore, let $I_\Gamma: \mathcal{C}^0(\Gamma) \rightarrow \mathbb{V}_\Gamma$ and $I_\Omega: \mathcal{C}^0(\Omega) \rightarrow \mathbb{V}_\Omega$ be the Lagrangian interpolant operators. The spatially discrete formulation is finally given by: find $B,Q: \mathbb{V}_\Omega \times [0,T]$ and $\Lambda,\Theta : \mathbb{V}_\Gamma \times [0,T] \rightarrow \mathbb{R}$ such that
\begin{equation}
\label{semidiscrete_formulation}
\begin{cases}
m_h^\Omega \left(\dot{B}, \Phi\right) + a_h^\Omega\left(B, \Phi\right) = m_h^\Omega \left(I_\Omega\ftilde_1(B), \Phi\right) - \psi_\eta m_h^\Gamma\left(I_\Gamma\ftilde_3(B, \Lambda, \Theta), \Phi\right);\\
m_h^\Omega \left(\dot{Q}, \Phi\right) + d_\Omega a_h^\Omega\left(Q, \Phi \right) = m_h^\Omega \left(I_\Omega \ftilde_2(Q),  \Phi\right) - d_\Omega\psi_\theta m_h^\Gamma\left(I_\Gamma\ftilde_4(Q, \Lambda, \Theta), \Phi\right);\\
m_h^\Gamma\left(\dot{\Lambda}, \Psi\right) + a_h^\Gamma\left(\Lambda, \Psi\right) = m_h^\Gamma\left(I_\Gamma\ftilde_3(B,\Lambda, \Theta), \Psi\right);\\
m_h^\Gamma \left( \dot{\Theta}, \Psi\right) + d_\Gamma a_h^\Gamma\left( \Theta,   \Psi\right) = m_h^\Gamma\left(I_\Gamma\ftilde_4(Q,\Lambda, \Theta), \Psi\right),
\end{cases}
\end{equation}
for all $\Phi : \mathbb{V}_\Omega \times [0,T] \rightarrow\mathbb{R}$ and $\Psi : \mathbb{V}_\Gamma \times [0,T] \rightarrow\mathbb{R}$.  If $N_\Gamma := \dim \mathbb{V}_\Gamma$ and $N_\Omega := \dim \mathbb{V}_\Omega$, let $\{\psi_i\}_{i=1}^{N_\Gamma}$ and $\{\varphi_i\}_{i=1}^{N_\Omega}$ be the Lagrangian bases of $\mathbb{V}_\Gamma$ and $\mathbb{V}_\Omega$, respectively.  We express the numerical solution $(B,Q,\Lambda,\Theta)$ in the Lagrange bases:
\begin{align}
\label{expansion_b}
&B(\bold x,t) = \sum_{i=1}^{N_\Omega} b_i(t) \varphi_i(\bold x), \qquad (\bold x,t) \in \Omega \times [0,T];\\
&Q(\bold x,t) = \sum_{i=1}^{N_\Omega} q_i(t) \varphi_i(\bold x), \qquad (\bold x,t) \in \Omega \times [0,T];\\
&\Lambda(\bold x,t) = \sum_{i=1}^{N_\Gamma} \lambda_i(t) \psi_i(\bold x), \qquad (\bold x,t) \in \Gamma \times [0,T];\\
\label{expansion_theta}
&\Theta(\bold x,t) = \sum_{i=1}^{N_\Gamma} \theta_i(t) \psi_i(\bold x), \qquad (\bold x,t) \in \Gamma \times [0,T],
\end{align}
where $b_i(t)$, $q_i(t)$, $\lambda_i(t)$, $\theta_i(t)$ are unknown time-dependent coefficients, which are collected in column vectors $\bold b(t), \bold q(t) \in \mathbb{R}^{N_\Omega}$, $\boldsymbol \eta(t), \boldsymbol \theta(t) \in \mathbb{R}^{N_\Gamma}$.  Following \cite{bsvem_parabolic}, we substitute \eqref{expansion_b}-\eqref{expansion_theta} into the spatially discrete formulation \eqref{semidiscrete_formulation},  and we obtain the following ODE system in vector form:
\begin{equation}
\label{spatially_discrete_vector_form}
\begin{cases}
M_\Omega \dot{\bold b} + A_\Omega \bold b = M_\Omega \ftilde_1(\bold b) - \psi_\eta R M_\Gamma \ftilde_3(\bold b, \boldsymbol \eta, \boldsymbol \theta);\\
M_\Omega \dot{\bold q} + d_\Omega A_\Omega \bold q = M_\Omega \ftilde_2(\bold q) - d_\Omega\psi_\theta R M_\Gamma \ftilde_4(\bold q, \boldsymbol \eta, \boldsymbol \theta);\\
M_\Gamma \dot{\boldsymbol \eta} + A_\Gamma \boldsymbol \eta = M_\Gamma \ftilde_3(\bold b, \boldsymbol \eta, \boldsymbol \theta);\\
M_\Gamma \dot{\boldsymbol \theta} + d_\Gamma A_\Gamma \boldsymbol \theta = M_\Gamma \ftilde_4(\bold q, \boldsymbol \eta, \boldsymbol \theta),
\end{cases}
\end{equation}
where the stiffness matrices $A_\Omega \in\mathbb{R}^{N_\Omega \times N_\Omega}$, $A_\Gamma \in\mathbb{R}^{N_\Gamma \times N_\Gamma}$, the lumped mass matrices $M_\Omega \in\mathbb{R}^{N_\Omega \times N_\Omega}$, $M_\Gamma\in\mathbb{R}^{N_\Gamma \times N_\Gamma}$ and the reduction matrix $R\in\mathbb{R}^{N_\Omega \times N_\Gamma}$ are defined as follows:
\begin{align}
&(A_\Omega)_{ij} := a_h^\Omega(\varphi_i, \varphi_j) \qquad (M_\Omega)_{ij} := m_h^\Omega (\varphi_i, \varphi_j), \qquad i,j=1,\dots, N_\Omega;\\
&(A_\Gamma)_{ij} := a_h^\Gamma(\psi_i, \psi_j), \qquad (M_\Gamma)_{ij} := m_h^\Gamma(\psi_i, \psi_j), \qquad i,j=1,\dots, N_\Gamma;\\
&R = \begin{bmatrix}
I_{N_\Gamma}\\
0
\end{bmatrix},
\end{align}
where $I_{N_\Gamma}$ is the identity of dimension $N_\Gamma$.

\subsection{Step 4: Fully discrete formulation in vector form}
Following \cite{bsvem_parabolic}, we discretize in time the ODE system \eqref{spatially_discrete_vector_form} with the IMEX Euler scheme.  Let $\tau > 0$ be the timestep and let $N_T = \left\lceil \frac{T}{\tau}\right\rceil$ be the number of timesteps. For all $n = 0, \dots, N_T - 1$,  the fully discrete solution $(\bold b^{(n)}, \bold q^{(n)}, \boldsymbol \eta^{(n)}, \boldsymbol \theta^{(n)})$ is found as follows:
\begin{equation}
\label{fully_discrete_vector}
\begin{cases}
(M_\Omega + \tau A_\Omega) \bold b^{(n+1)} =  M_\Omega\bold b^{(n)} + \tau \left(M_\Omega\bold f_1^{(n)} - \psi_\eta R M_\Gamma \bold f_3^{(n)}\right);\\
(M_\Omega + d_\Omega \tau A_\Omega) \bold q^{(n+1)} =  M_\Omega\bold q^{(n)} + \tau \left(M_\Omega\bold f_2^{(n)} - \psi_\theta d_\Omega R M_\Gamma \bold f_4^{(n)}\right);\\
(M_\Gamma + \tau A_\Gamma) \boldsymbol \eta^{(n+1)} = M_\Gamma\boldsymbol \eta^{(n)} + \tau M_\Gamma\bold f_3^{(n)};\\
(M_\Gamma + d_\Gamma \tau A_\Gamma) \boldsymbol \theta^{(n+1)} = M_\Gamma\boldsymbol \theta^{(n)} + \tau M_\Gamma\bold f_4^{(n)},\\
\end{cases}
\end{equation}
where
\begin{align*}
&\bold f_1^{(n)} := \ftilde_1(\bold b^{(n)});\\
&\bold f_2^{(n)} := \ftilde_2(\bold q^{(n)});\\
&\bold f_3^{(n)} := \ftilde_3(\bold b^{(n)}, \boldsymbol \eta^{(n)}, \boldsymbol \theta^{(n)});\\
&\bold f_4^{(n)} := \ftilde_4(\bold q^{(n)}, \boldsymbol \eta^{(n)}, \boldsymbol \theta^{(n)}).
\end{align*}
The fully discrete formulation \eqref{fully_discrete_vector} is composed of four linear algebraic systems that can be solved independently of each other at each time step. Of these four linear systems, two have dimension $N_\Omega$, while the other two have dimension $N_\Gamma$.  If the cube $\Omega$ is discretised with a Cartesian mesh with $N_x \in\mathbb{N}$ gridpoints along each dimension, then $N_\Omega = N_x^3$, which makes the linear systems in \eqref{fully_discrete_vector} computationally prohibitive to solve.  An extremely efficient approach to address this issue is the so-called Matrix-Oriented Finite Element Method (MOFEM) \cite{mofem3d}, which exploits the Cartesian structure of the grid to translate the linear systems in \eqref{fully_discrete_vector} into tensor equations of much lower size. We will show a numerical solution to our problem carried out with MOFEM  in Section \ref{sec:simulations_turing}. However, we will show that, given the particular nature of the considered PDE problem, an even more efficient solver is given by the BSVEM on a bespoke mesh.

\subsection{Bespoke BS-VEM mesh for the BS-DIB model}
\label{sec:bespoke_mesh}
\begin{figure}[h!]
\begin{center}
\subfigure[Nodes of the adaptive mesh]{\includegraphics[scale=.45]{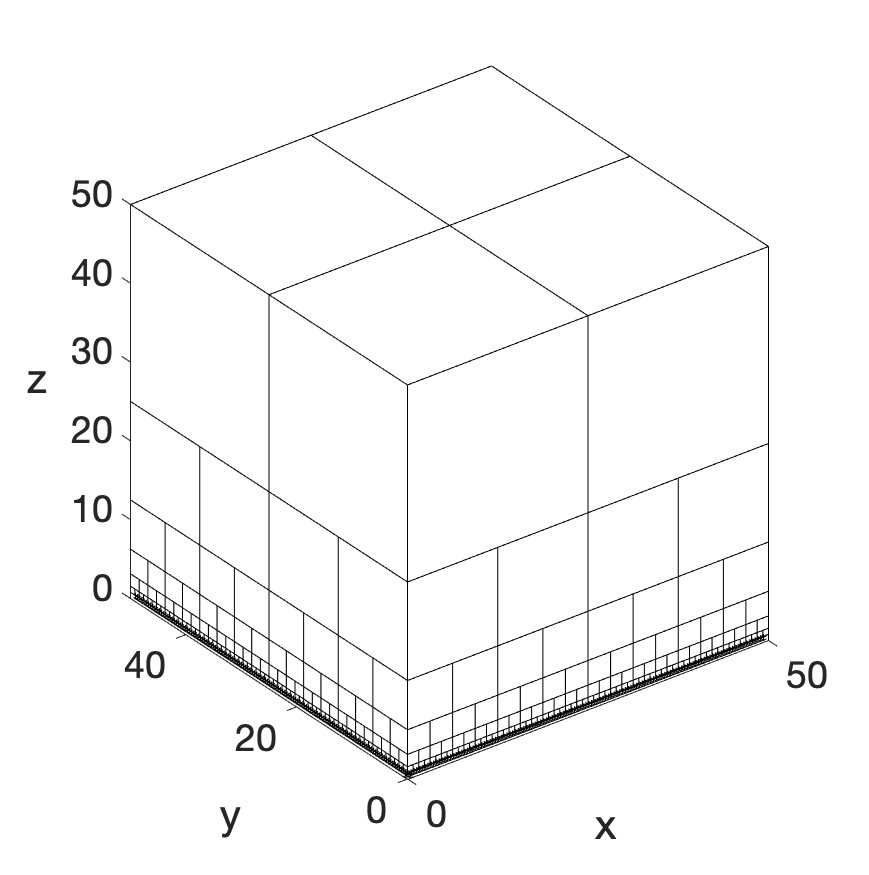}}
\subfigure[Cubic element with hanging nodes]{\includegraphics[scale=.45]{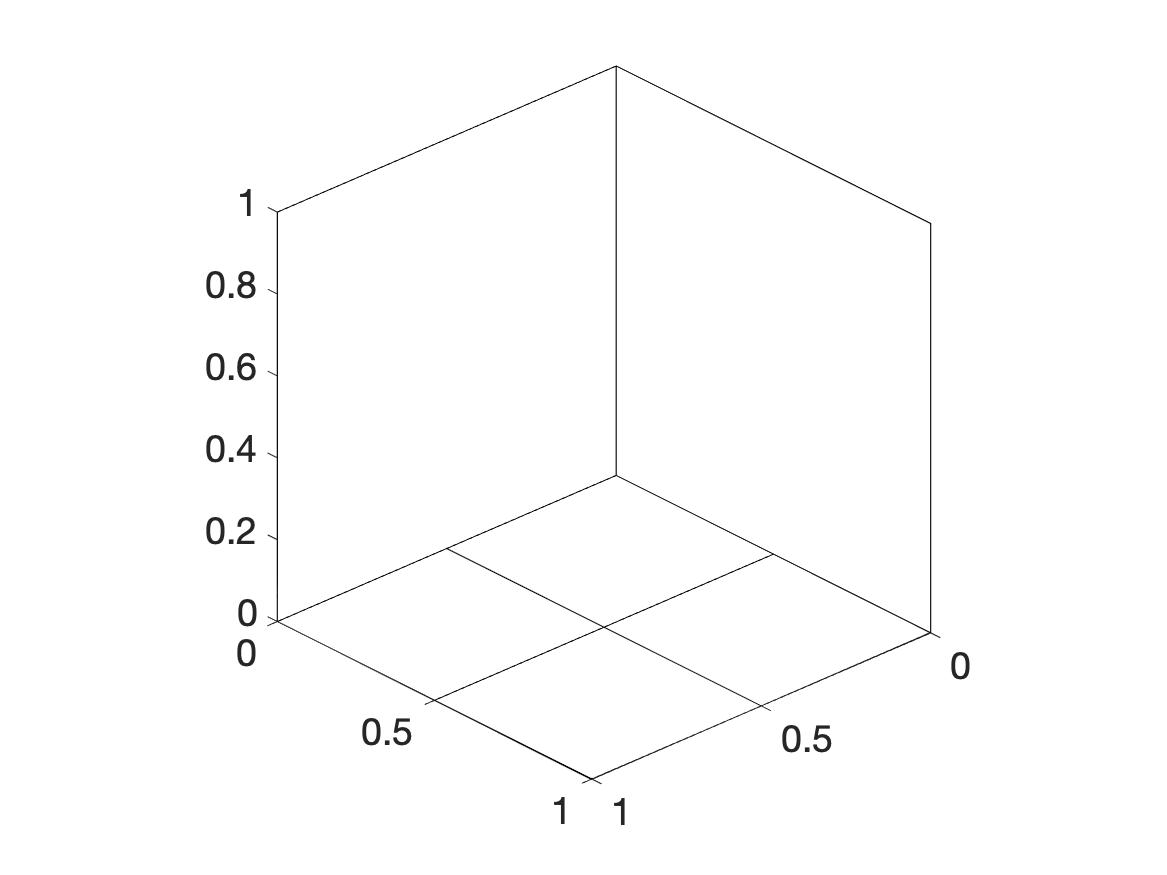}}
\end{center}
\caption{Graded polyhedral mesh used in the BS-VEM approximation of the model \eqref{model}.}
\label{fig:graded_mesh}
\end{figure}

\begin{figure}[h!]
\begin{center}
\subfigure[Comparison between the BS-DIB model \eqref{model} and the DIB model \eqref{model2d}: surface component $\eta$ at the final time $T=50$. ]{\includegraphics[scale=.4]{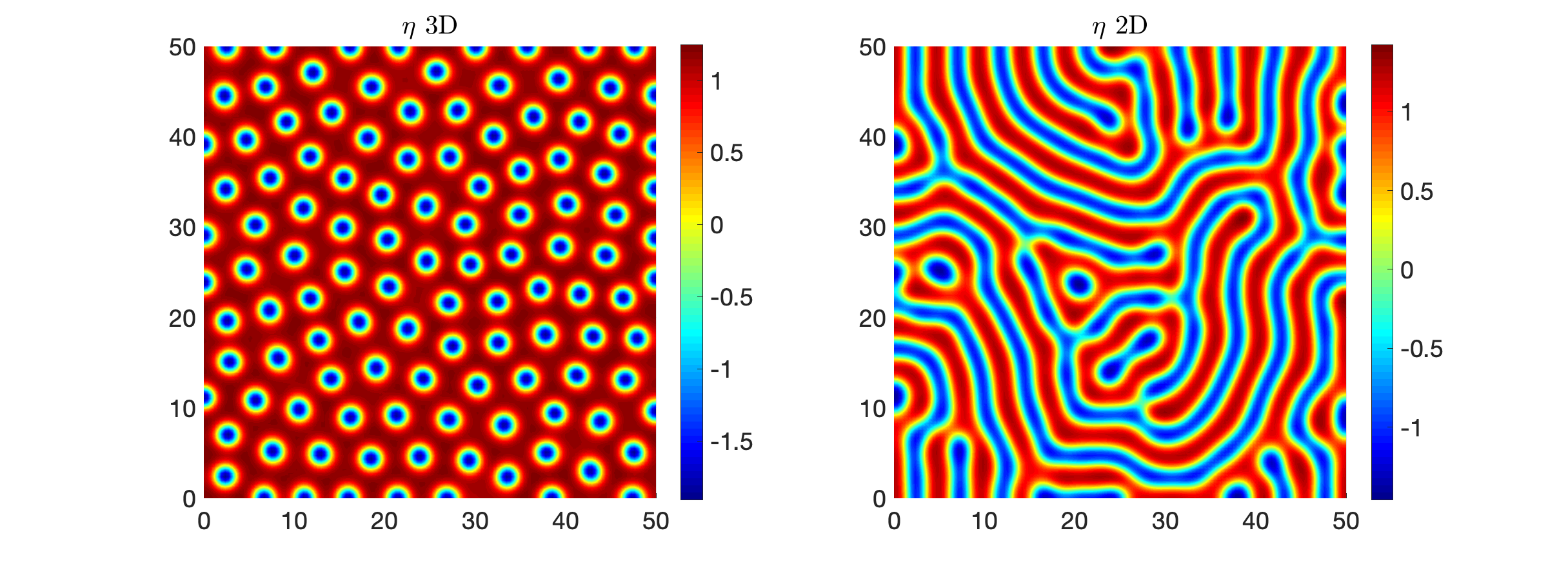}\label{fig:experimentD3_MOFEM_comparison_eta}}
\subfigure[Coupled 3D DIB model \eqref{model}: bulk component $b$ and surface component $\eta$ at the final time $T=50$. In the bulk, the spatial pattern is confined to a neighborhood of the surface $\Gamma$. ]{\includegraphics[scale=.4]{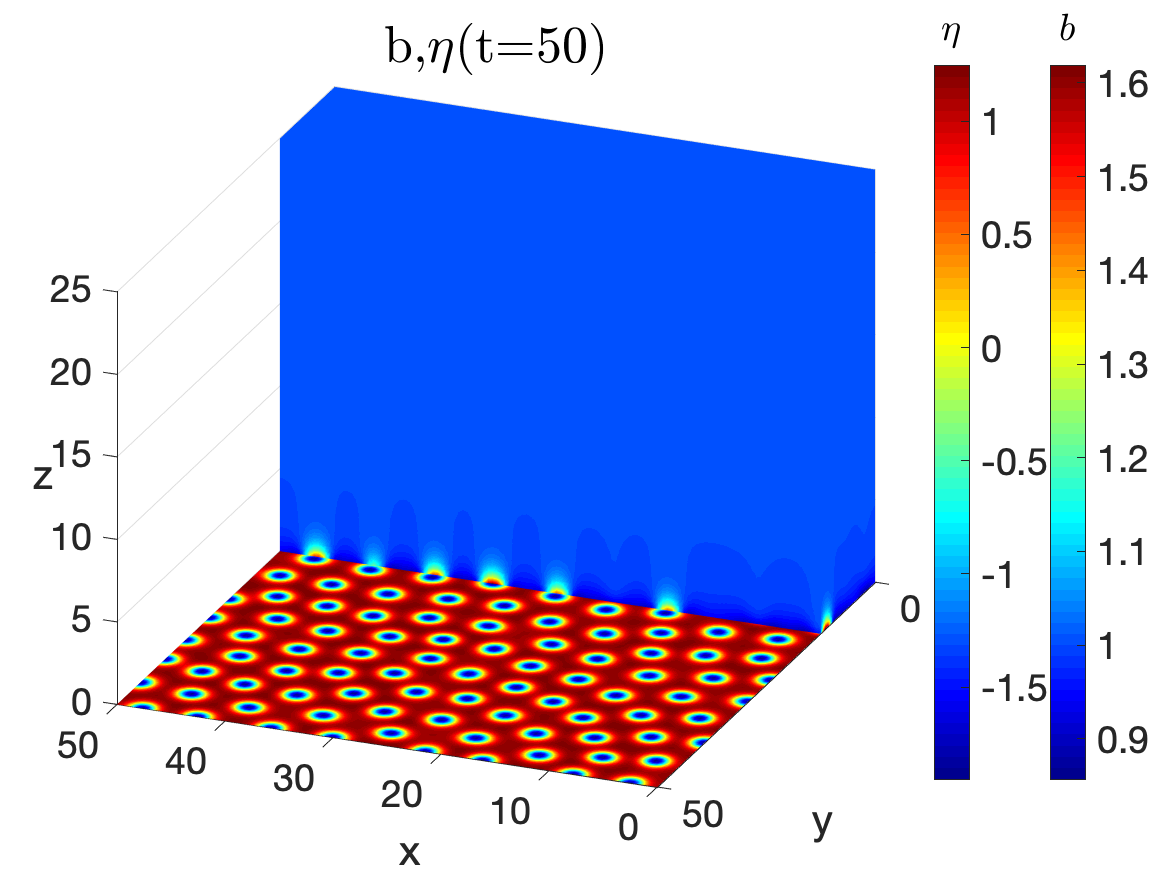}\label{fig:experimentD3_MOFEM_bs}}
\subfigure[Coupled 3D DIB model \eqref{model}: increment of surface component $\eta$ over time.]{\includegraphics[scale=.4]{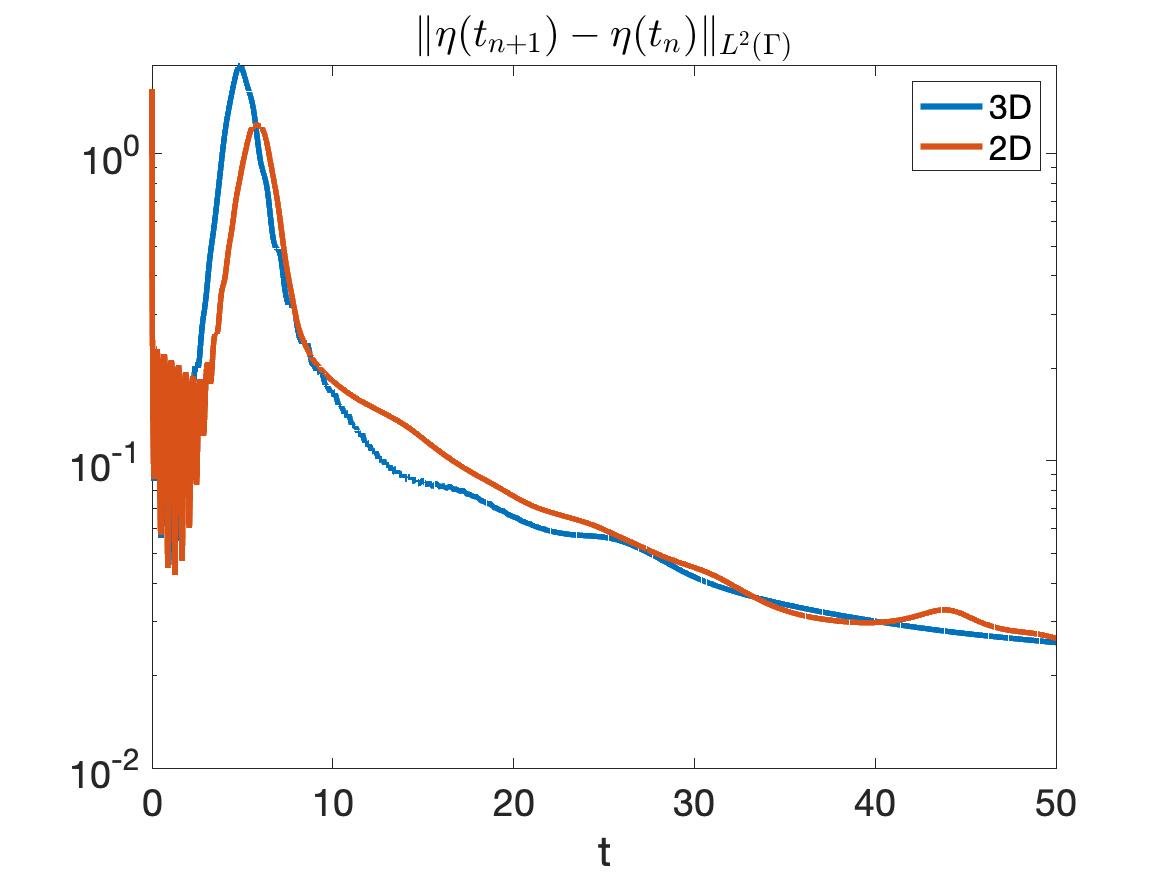}\label{fig:experimentD3_MOFEM_increment}}
\caption{Simulation obtained with the MO-FEM approach. Bifurcation and coupling  parameters as in Table \ref{tab:experiments_recap} (Experiment D3). The BS-DIB model \eqref{model} shows a reversed spots pattern, while the DIB model \eqref{model2d} exhibits a labyrinth pattern.  In the BS-DIB model \eqref{model}, the bulk components $(b,q)$ exhibit a spatial pattern only in a neighborhood of the surface $\Gamma$. This suggests the usage of a graded mesh that is highly refined close to $\Gamma$ and much coarser away from $\Gamma$.}
\label{fig:experimentD3_MOFEM}
\end{center}
\end{figure}

Since the domain of the model problem \eqref{model} is a cube, then it would be natural to choose an efficient numerical method that exploits the structure of Cartesian grids, such as the Matrix-Oriented Finite Element Method (MOFEM) \cite{mofem3d}.  By choosing the model parameters and the timestep as in Table \ref{tab:experiments_recap}, Experiment D1, and a Cartesian grid of $128\times 128 \times 128 \approx 2.1$e+6 nodes, MOFEM produces the numerical solution shown in Fig.  \ref{fig:experimentD3_MOFEM}.  In keeping with the physical meaning of the parameter choice, we can observe that the bulk components $(b,q)$ exhibit spatial patterns only in the proximity of the surface $\Gamma$, and become approximately constant away from $\Gamma$. This suggests that a uniform Cartesian cubic grid would be unnecessary fine away from $\Gamma$. Hence, we apply the BS-VEM with a graded cubic mesh that is highly refined close to $\Gamma$ and gradually becomes coarser as the distance from $\Gamma$ increases. Such graded polyhedral mesh is depicted in Fig.  \ref{fig:graded_mesh}(a). This grid is composed of two layers of $128\times 128$ cubic elements close to $\Gamma$, and five layers of gradually larger ``false cubes''. Such a false cube is depicted Fig.  \ref{fig:graded_mesh}(b) and is actually an \emph{ennahedron}: a polyhedron with nine faces and thirteen vertices.  Specifically,  such ennahedron is a cube the bottom face of which is split into four smaller square faces (highlighted in purple in Fig. \ref{fig:graded_mesh}(b)). The proposed graded mesh provides a finer discretisation for the surface $\Gamma$ ($129\times 129$ nodes on $\Gamma$) than the Cartesian grid ($100\times 100$ nodes on $\Gamma$) used for the MOFEM and, at the same time has much less nodes (approximately $5,5$e+4 versus $1e$+6), resulting in a discrete problem of much smaller size and shorter computational times (approximately 89 minutes for MO-FEM and 39 minutes for BS-VEM).

\section{Numerical experiments}
\label{sec:simulations_turing}
We shall present seven numerical experiments to compare the DIB model \eqref{model2d} with the novel BSDIB model \eqref{model}.  The seven experiments differ from each other for appropriate choices of the model parameters.  The first four experiments, called T1 through T4,  show that for various choices of the bifurcation parameters, the BSDIB model shows Turing patterns while the DIB model does not. This seems to suggest that the BSDIB model has a larger Turing space than the DIB model.  As mentioned above, a theoretical analysis of the Turing space for the BSDIB model is outside the scope of this work. The latter three experiments, called D1 through D3,  show that when the DIB model exhibits Turing patterns, the BSDIB model still exhibits Turing patterns, but not necessarily of the same morphological class.  All the experiments are carried out on a cubic domain of edge length $L=50$ on the polyhedral mesh described in Section \ref{sec:bespoke_mesh}. The final time $T$ and the timestep $\tau$ also differ for each experiments according to the stiffness of the problem and the timescale of the dynamics. A recap of the numerical experiments and the respective parameters is given in Table \ref{tab:experiments_recap}.

\begin{table}[h!]
{\small
\begin{tabular}{c|c|c|c|c|c|c|c|c|c}
Experiment & $A_2$ & $B$ & $C$ & $\gamma$ & $\psi_\eta = \psi_\theta$ & $T$ & $\tau$ & Pattern 3D & Pattern 2D\\
\hline
T1 & 1 & 50 & 10 & 0 & 0.5 & 200 & 2e-3 & Thin worms & homogeneous\\
T2 & 1 & 75 & 5 & 0 & 0.3 & 200 & 2e-3 & Stripes & homogeneous\\
T3 & 1 & 35 & 15 & 0 & 0.5 & 100 & 5e-3 & Holes & homogeneous\\
T4 & 1 & 30 & 20 & 0 & 0.5 & 50 & 5e-3 & homogeneous (different from 2D) & homogeneous\\
D1 & 1 & 66 & 3 & 0.2 & 0.1 & 200 & 5e-3 & Holes and worms & Labyrinth\\
D2 & 1 & 66 & 3 & 0.2 & 0.15 & 50 & 2e-3 & Holes and small worms & Labyrinth\\
D3 & 1 & 66 & 3 & 0.2 & 0.2 & 50 & 2e-3 & Holes & Labyrinth\\
D4 & 1 & 3 & 30 & 0.2 & 0.1 & 200 & 5e-3 & Holes & Bigger Holes
\end{tabular}}
\caption{A recap of the numerical experiments.}
\label{tab:experiments_recap}
\end{table}

\begin{figure}[h!]
\begin{center}
\subfigure[Surface component $\eta$ at the final time $T=200$. ]{\includegraphics[scale=.4]{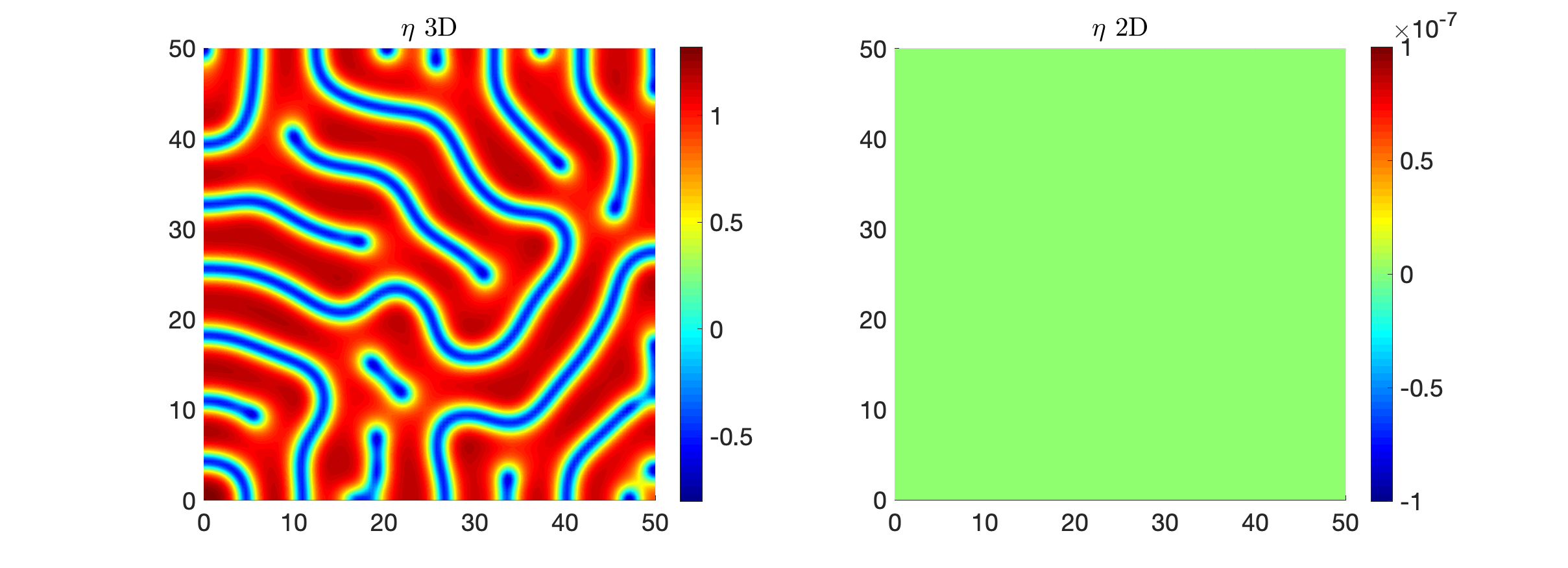}\label{fig:experimentT1_comparison_eta}}
\subfigure[Surface component $\theta$ at the final time $T=200$. ]{\includegraphics[scale=.4]{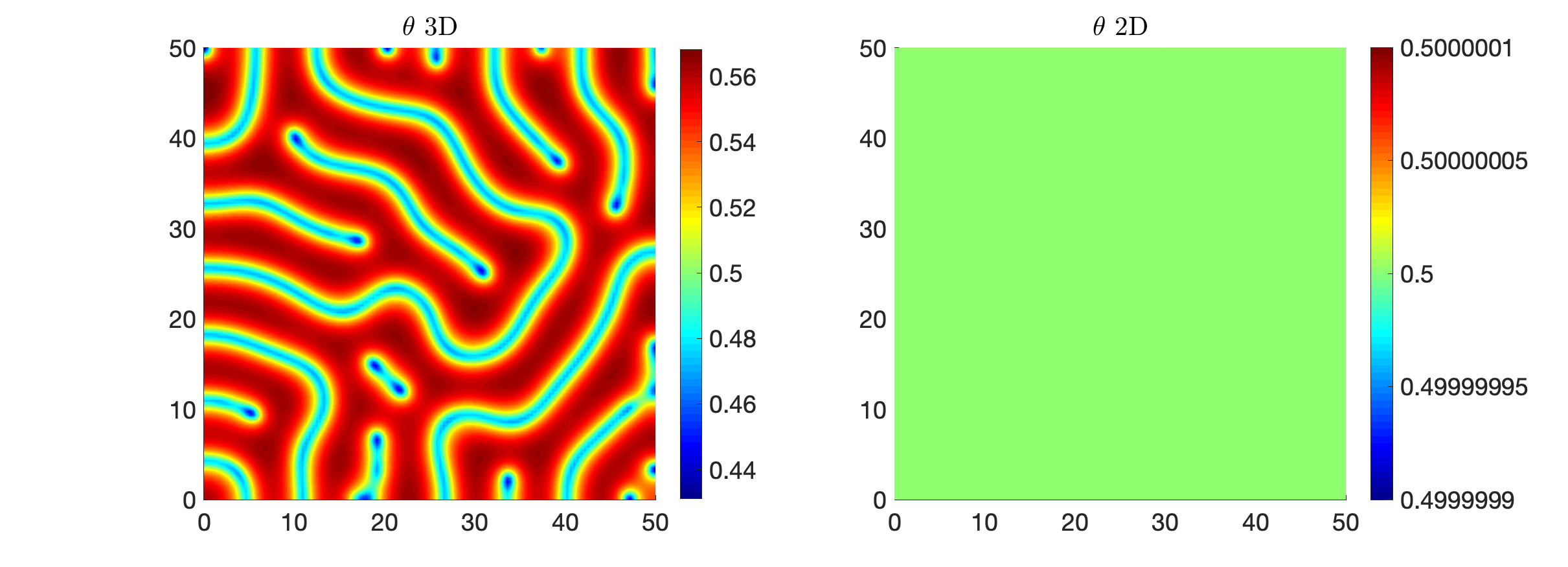}\label{fig:experimentT1_comparison_theta}}
\subfigure[BS-DIB model \eqref{model}: bulk component $b$ and surface component $\eta$ at final time $T=200$.  The $z$-axis is not in scale to help visualization.]{\includegraphics[scale=.4]{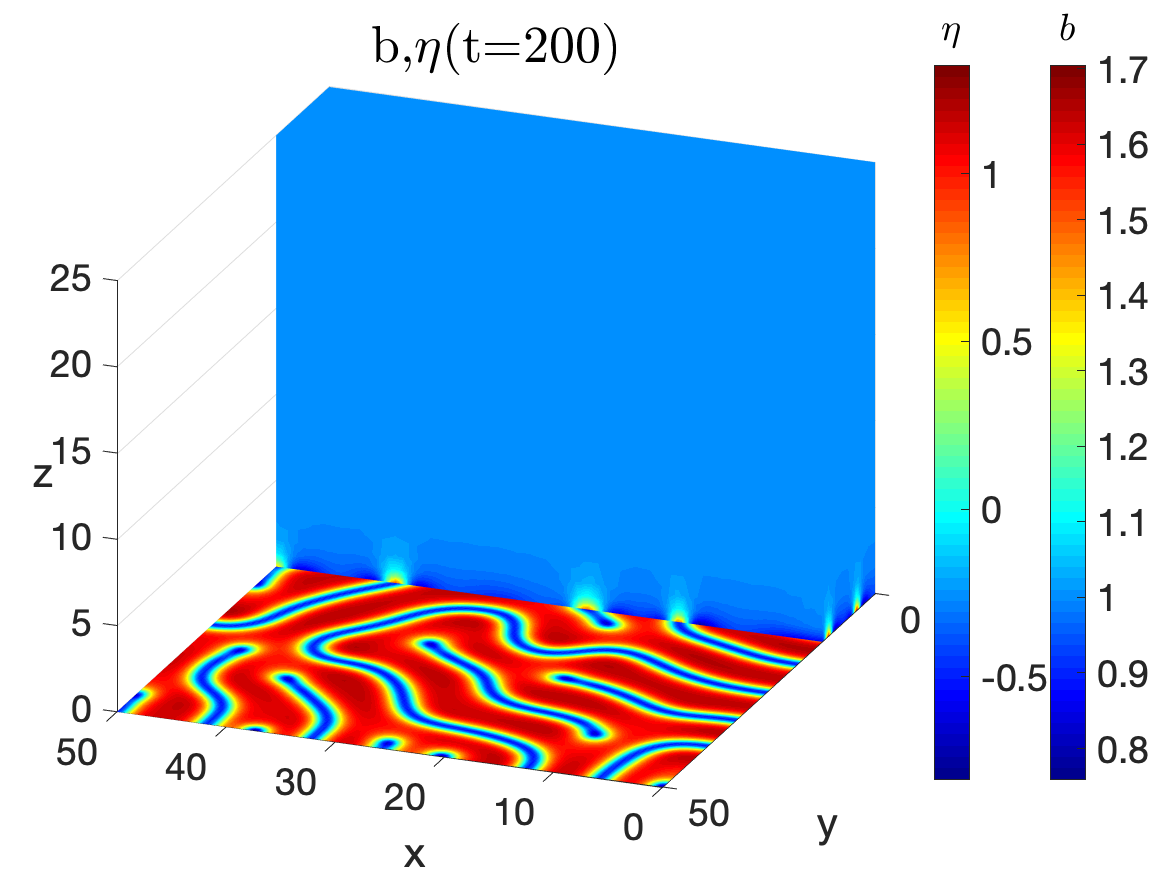}\label{fig:experimentT1_bs}}
\subfigure[Increment of surface component $\eta$ over time.]{\includegraphics[scale=.4]{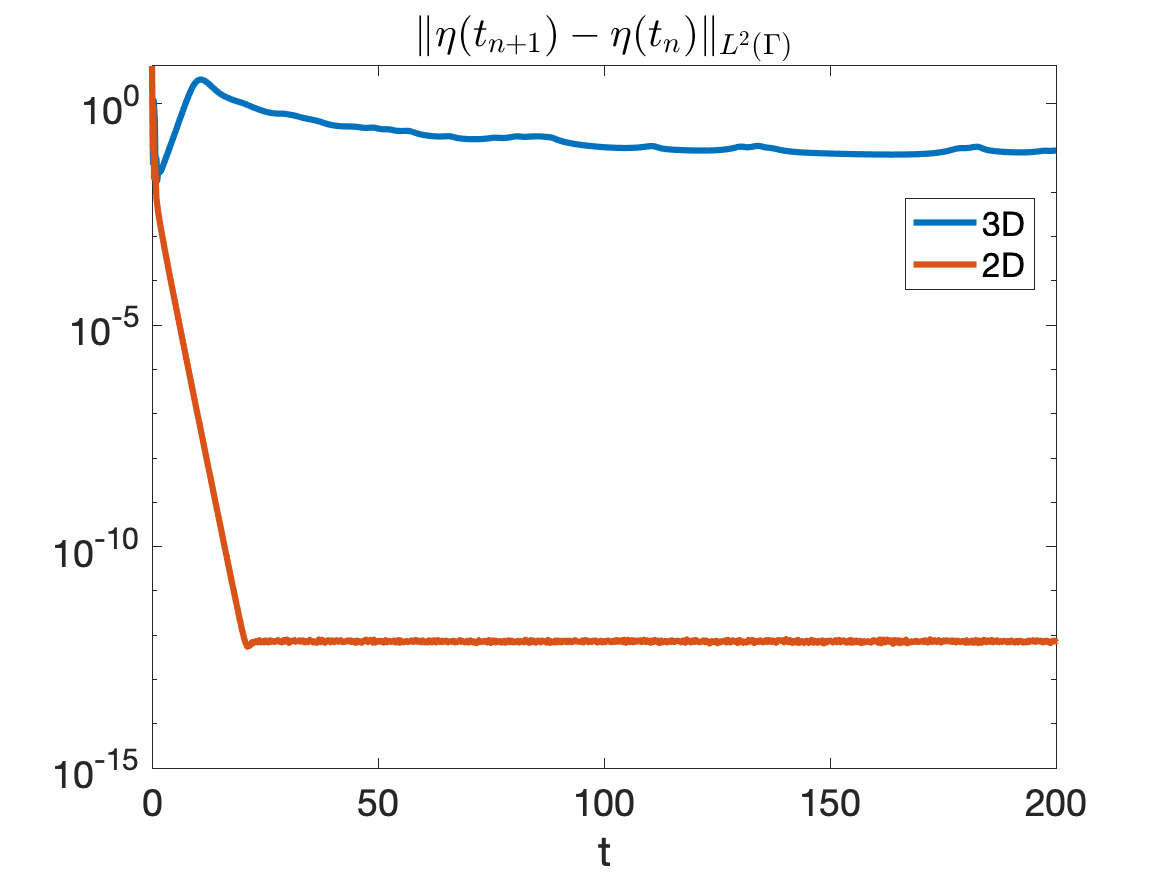}\label{fig:experimentT1_increment}}
\caption{Simulation T1.  Bifurcation and coupling  parameters as in Table \ref{tab:experiments_recap}. The BSDIB model \eqref{model} shows a worm pattern, while the DIB model \eqref{model2d} reaches the homogeneous steady state.}
\label{fig:experimentT1}
\end{center}
\end{figure}

\begin{figure}[h!]
\begin{center}
\subfigure[Surface component $\eta$ at the final time $T=200$. ]{\includegraphics[scale=.4]{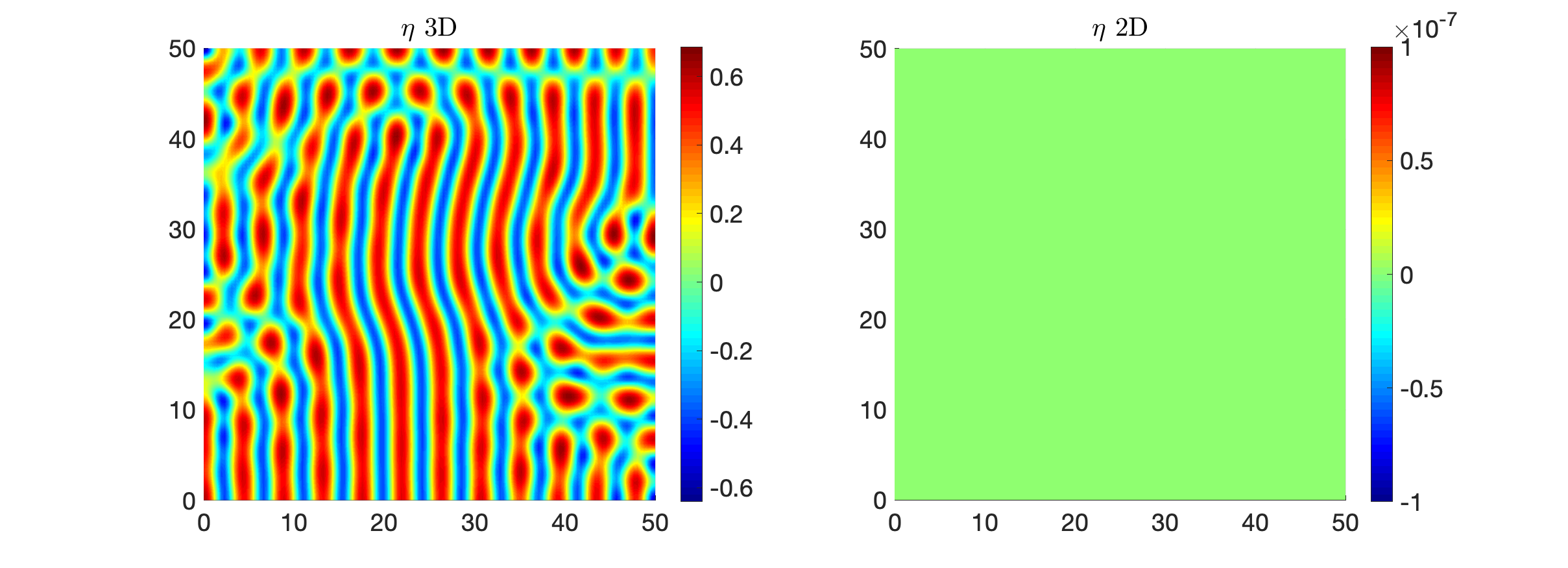}\label{fig:experimentT2_comparison_eta}}
\subfigure[Surface component $\eta$ at the final time $T=200$. ]{\includegraphics[scale=.4]{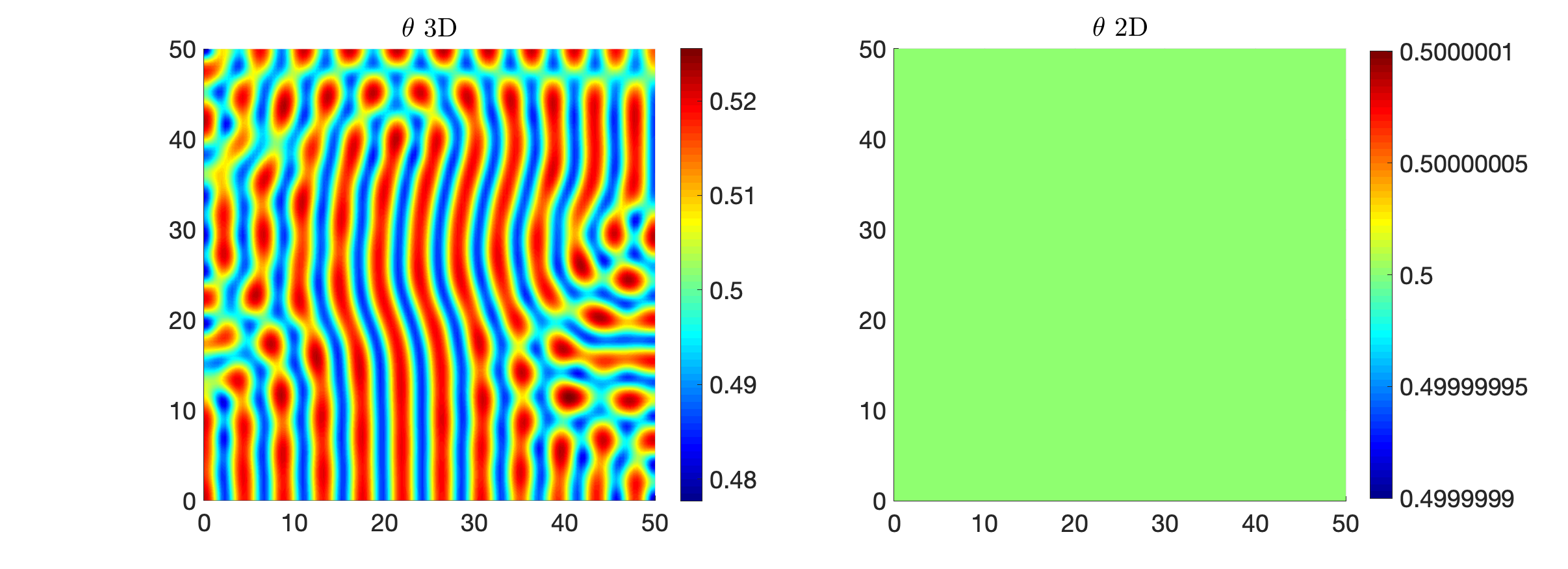}\label{fig:experimentT2_comparison_theta}}
\subfigure[BS-DIB model \eqref{model}: bulk component $b$ and surface component $\eta$ at the final time $T=200$. ]{\includegraphics[scale=.4]{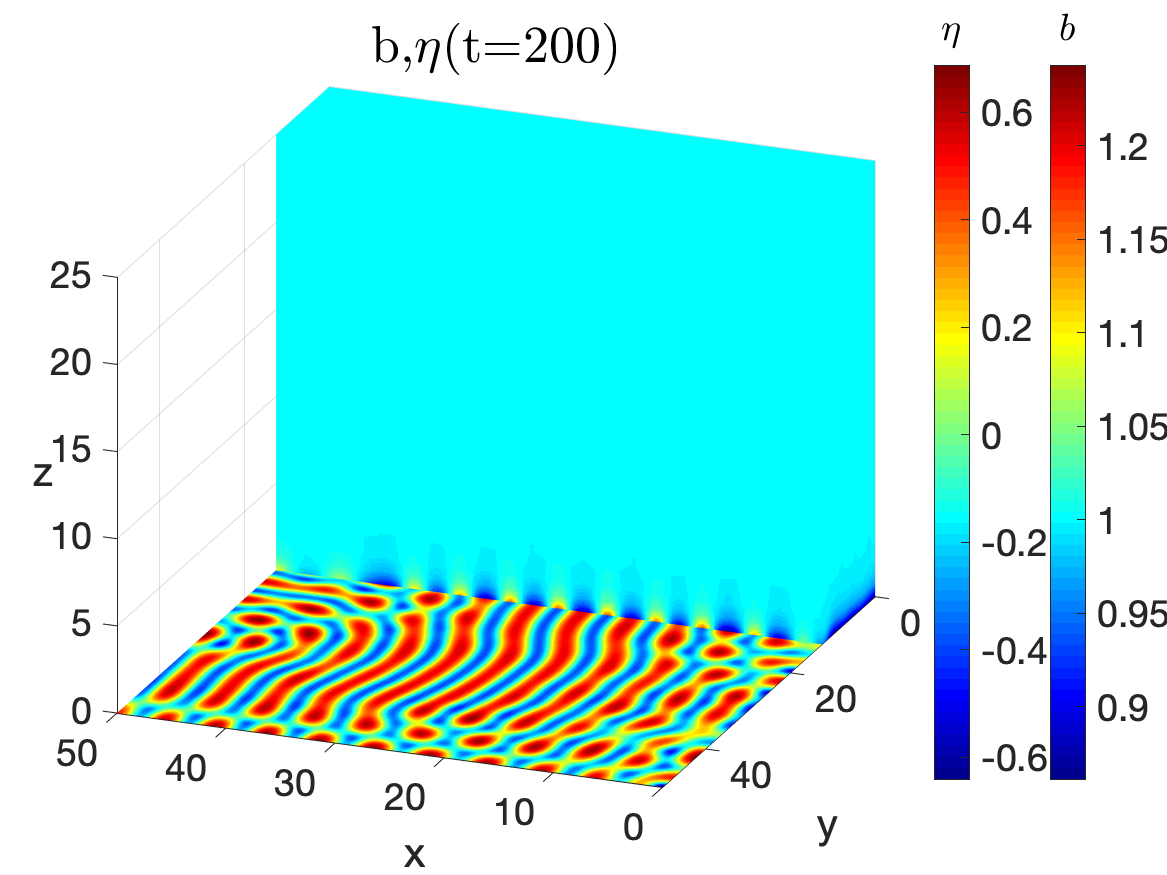}\label{fig:experimentT2_bs}}
\subfigure[Increment of surface component $\eta$ over time.]{\includegraphics[scale=.4]{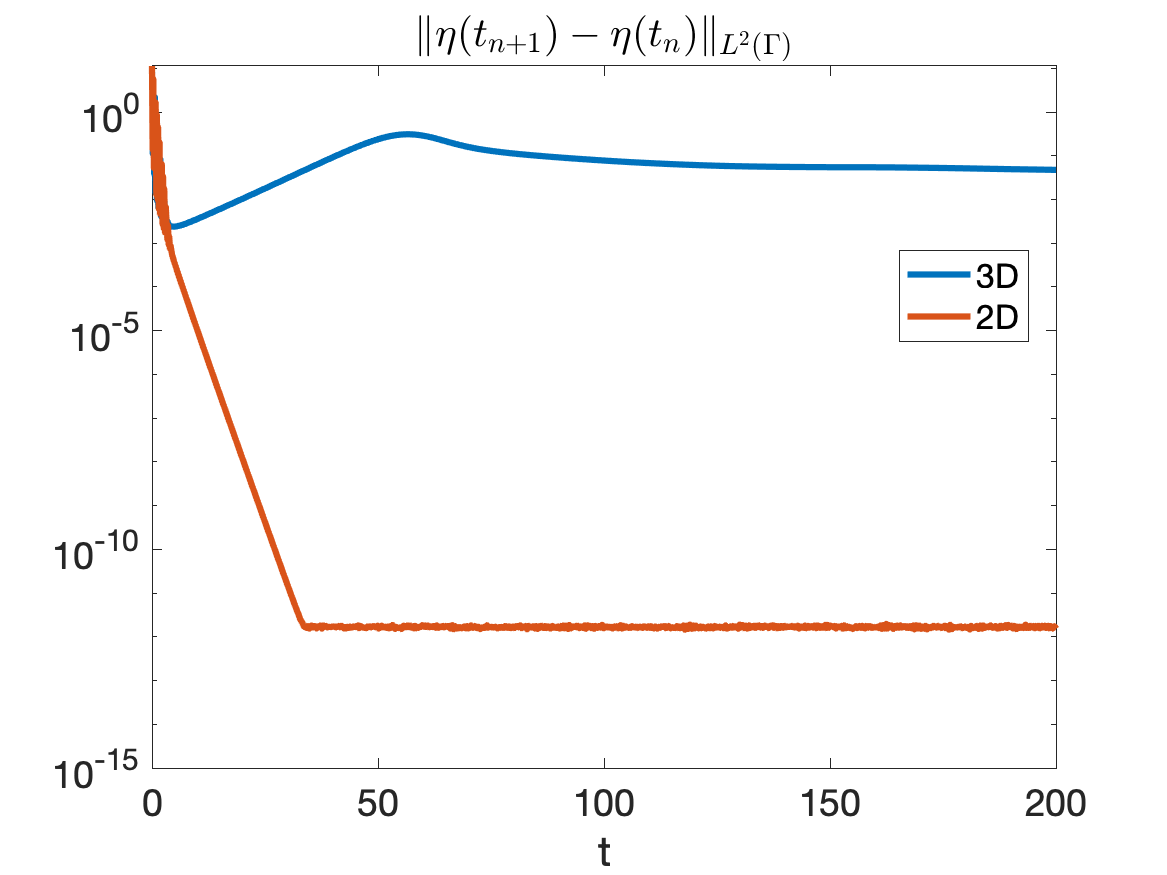}\label{fig:experimentT2_increment}}
\caption{Simulation T2.  Bifurcation and coupling  parameters as in Table \ref{tab:experiments_recap}. The BSDIB model \eqref{model} shows a slow-to-stabilize stripe pattern, while the DIB model \eqref{model2d} reaches the homogeneous steady state.}
\label{fig:experimentT2}
\end{center}
\end{figure}

\begin{figure}[h!]
\begin{center}
\subfigure[Surface component $\eta$ at the final time $T=200$. ]{\includegraphics[scale=.4]{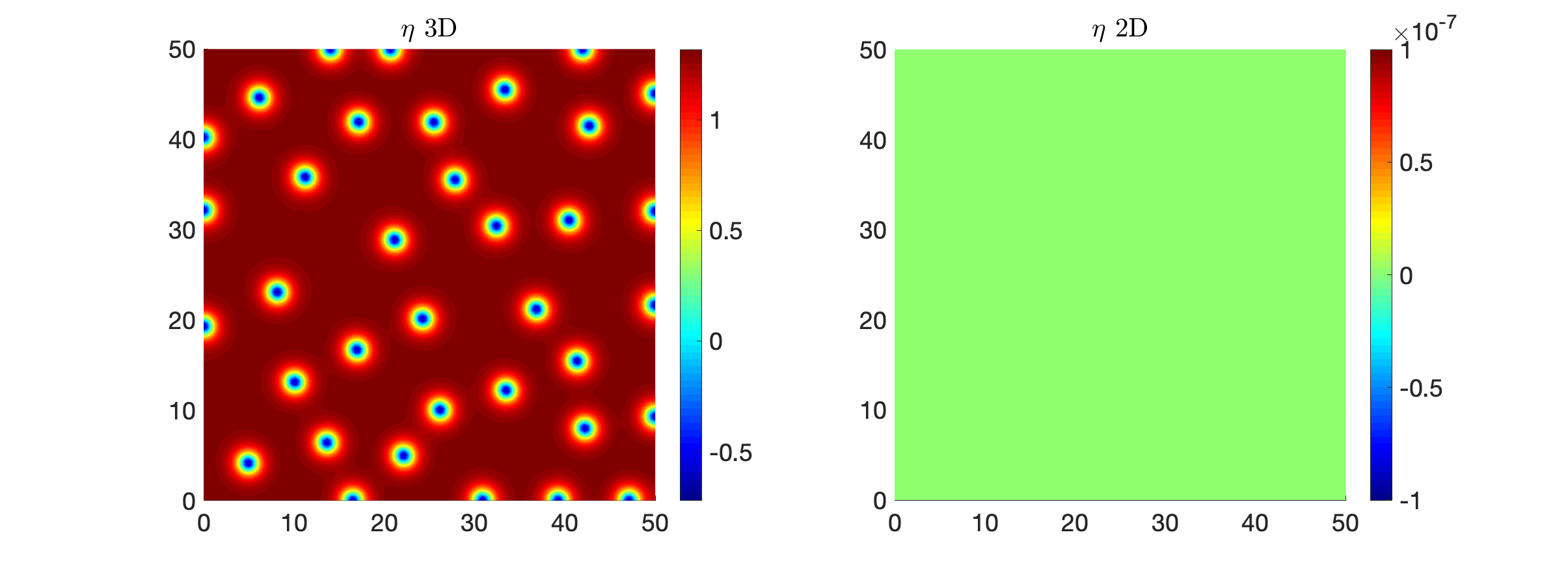}\label{fig:experimentT3_comparison_eta}}
\subfigure[Surface component $\eta$ at the final time $T=200$. ]{\includegraphics[scale=.4]{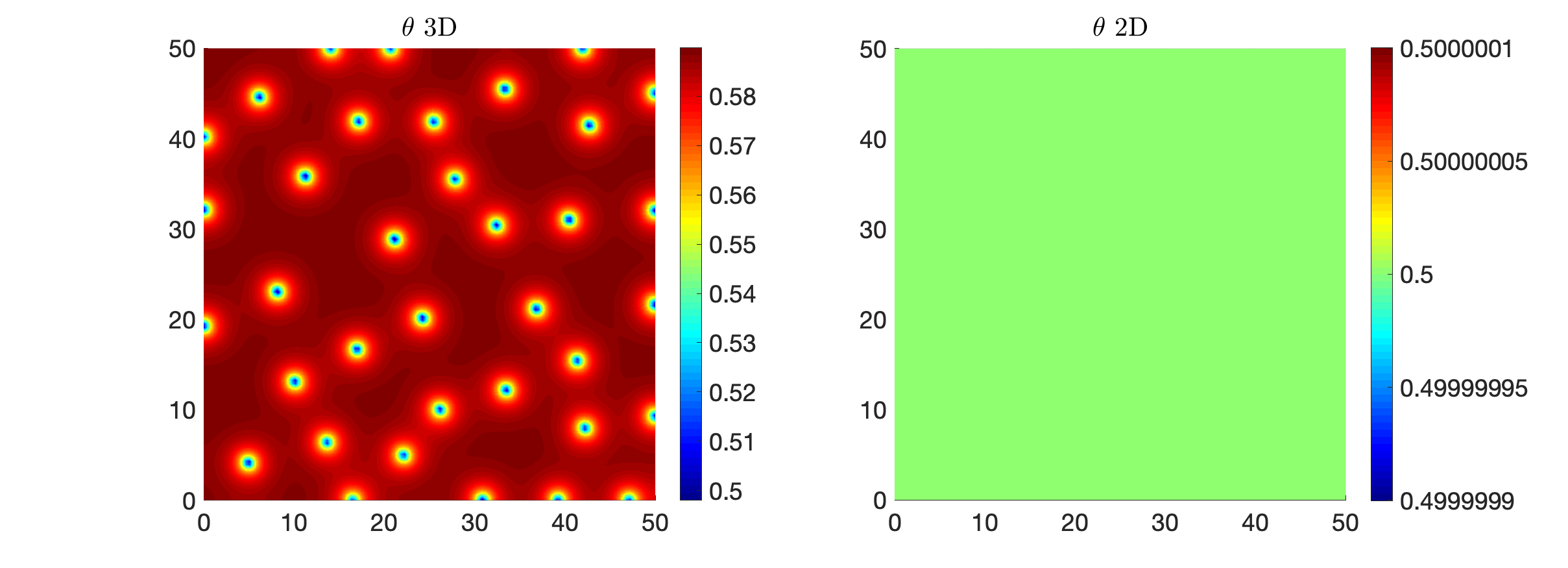}\label{fig:experimentT3_comparison_theta}}
\subfigure[BS-DIB model \eqref{model}: bulk component $b$ and surface component $\eta$ at the final time $T=200$. ]{\includegraphics[scale=.4]{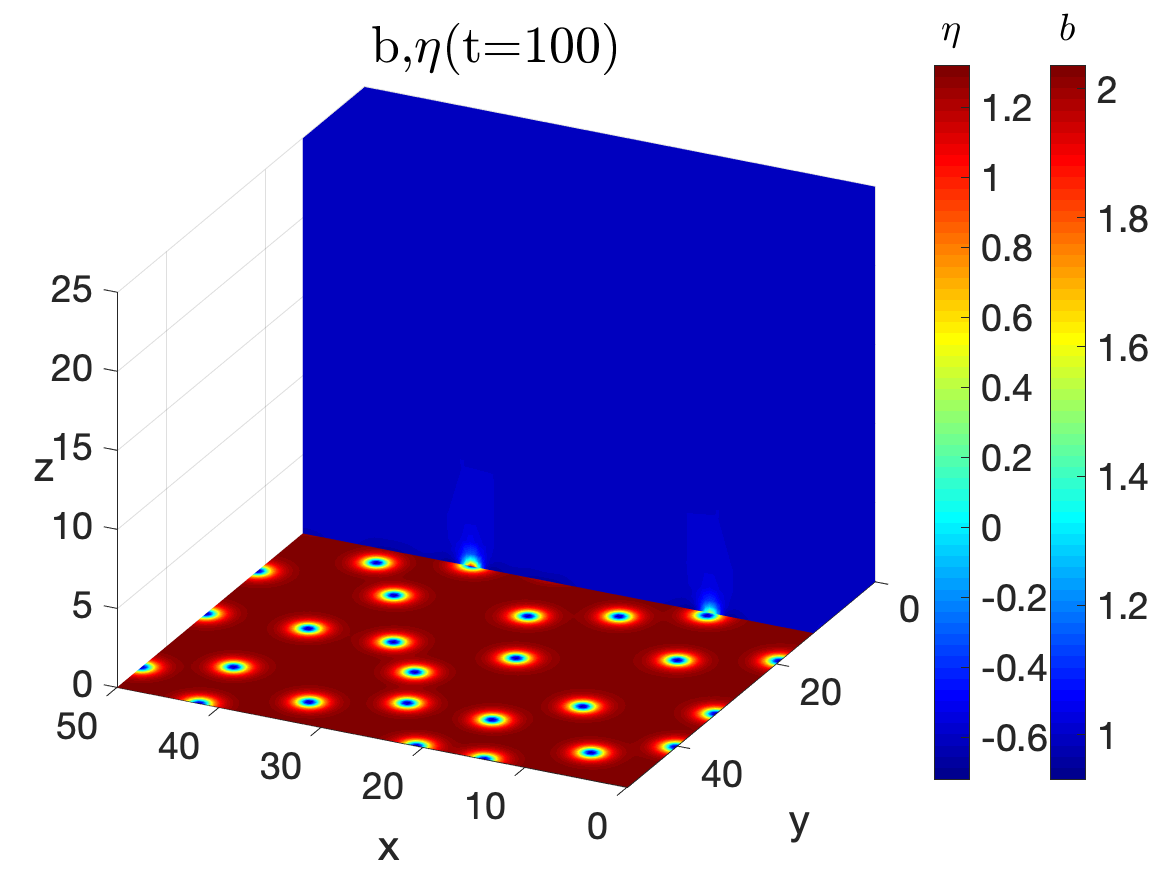}\label{fig:experimentT3_bs}}
\subfigure[Increment of surface component $\eta$ over time.]{\includegraphics[scale=.4]{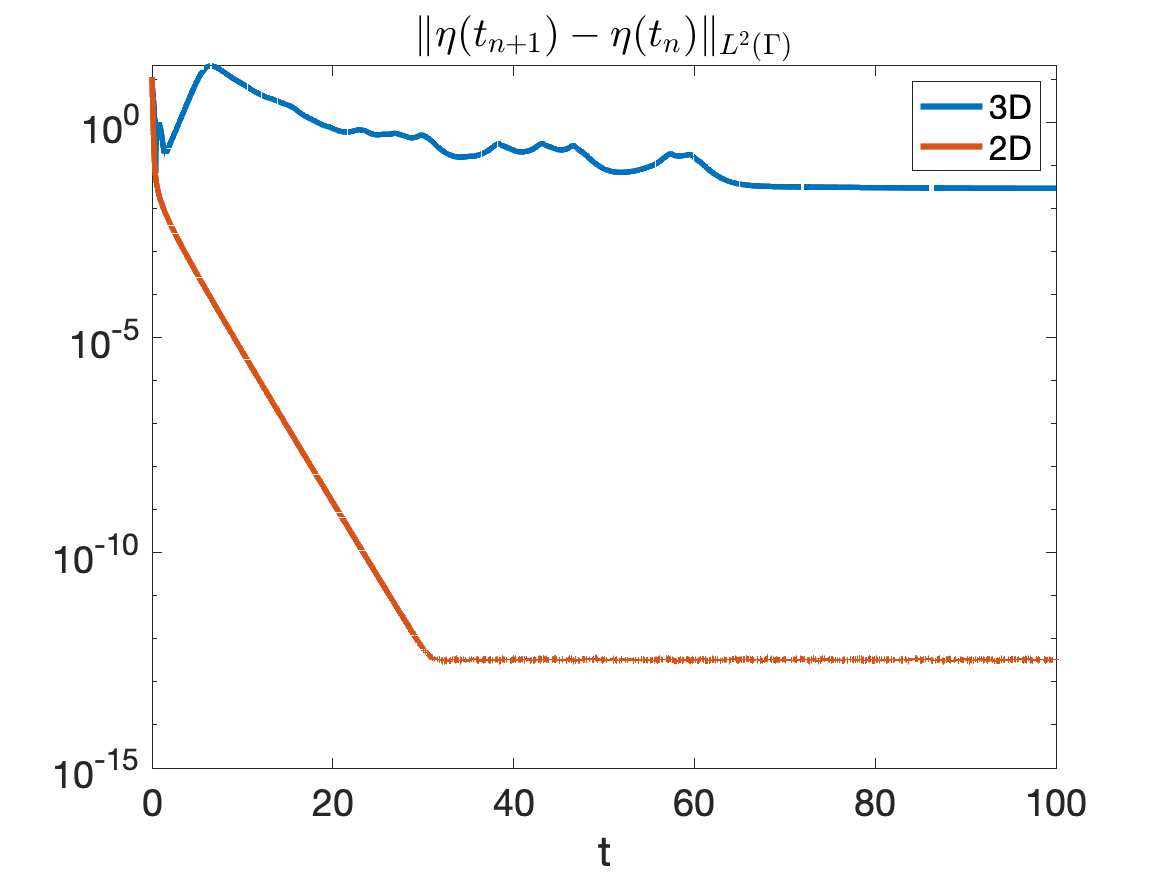}\label{fig:experimentT3_increment}}
\caption{Simulation T3.  Bifurcation and coupling  parameters as in Table \ref{tab:experiments_recap}. The BSDIB model \eqref{model} shows a reversed spots pattern, while the DIB model \eqref{model2d} reaches the homogeneous steady state.}
\label{fig:experimentT3}
\end{center}
\end{figure}

\begin{figure}[h!]
\begin{center}
\subfigure[Surface component $\eta$ at the final time $T=50$. ]{\includegraphics[scale=.4]{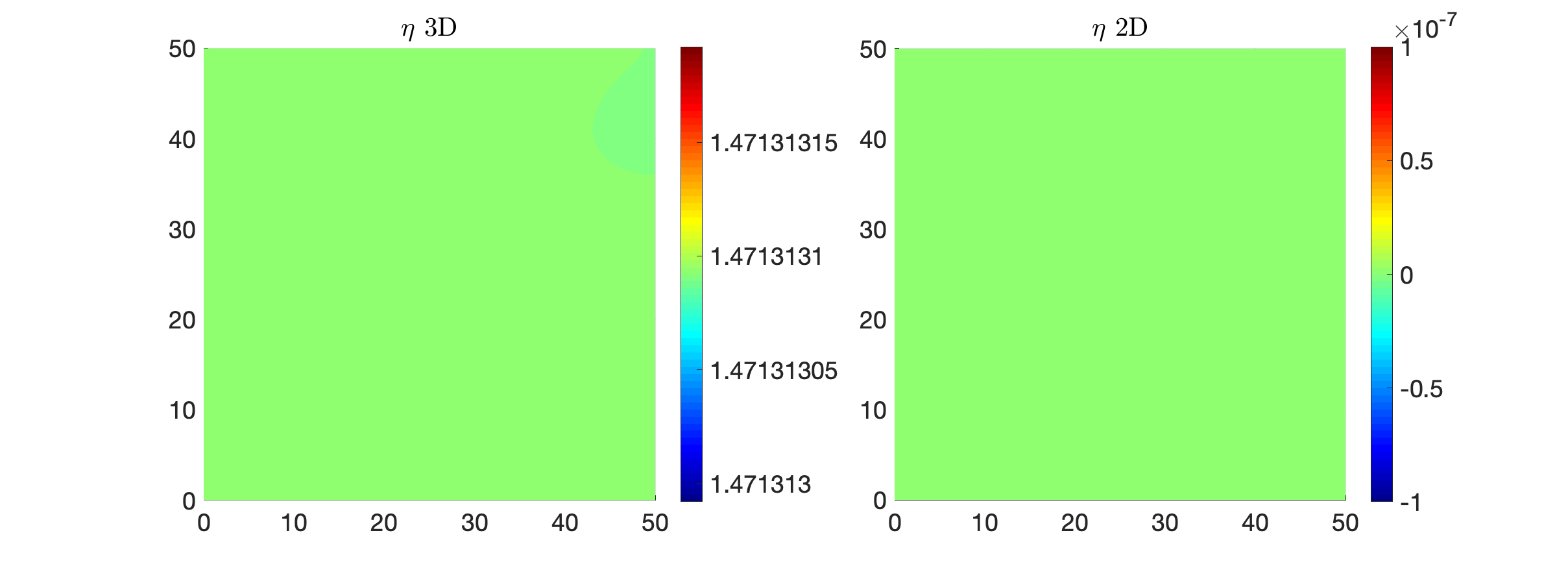}\label{fig:experimentT4_comparison_eta}}
\subfigure[Surface component $\eta$ at the final time $T=50$. ]{\includegraphics[scale=.4]{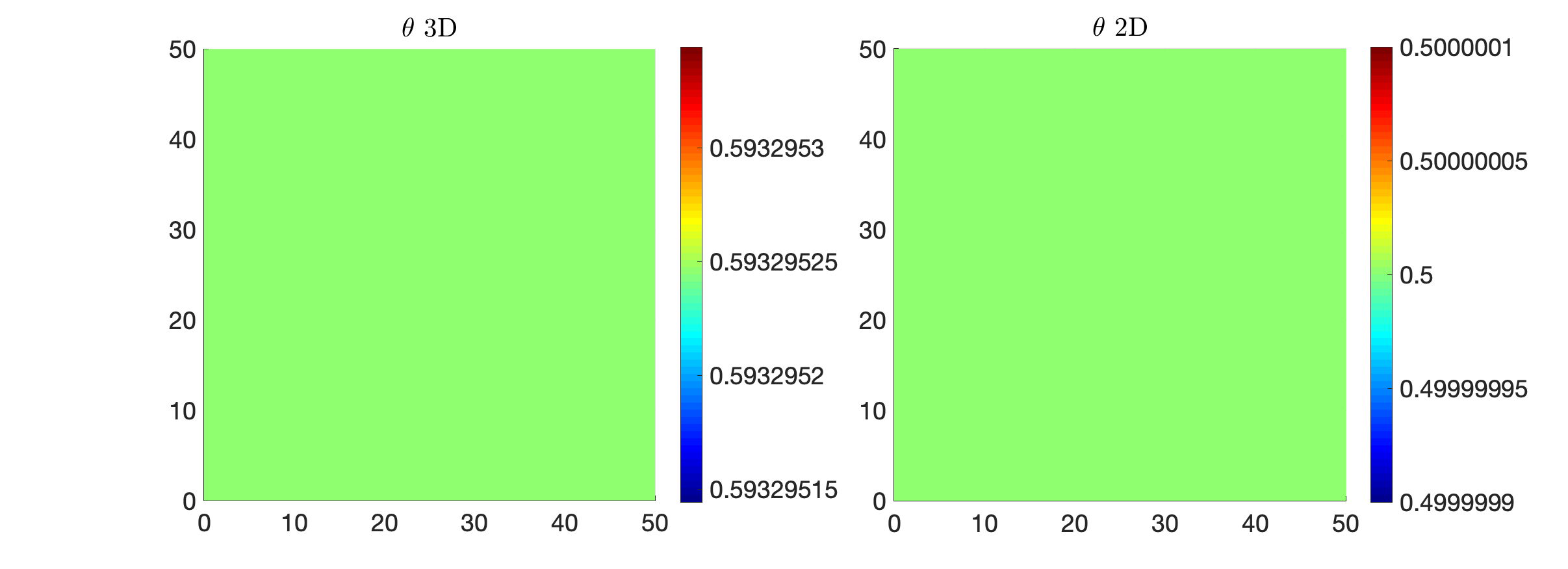}\label{fig:experimentT4_comparison_theta}}
\subfigure[Increment of surface component $\eta$ over time.]{\includegraphics[scale=.4]{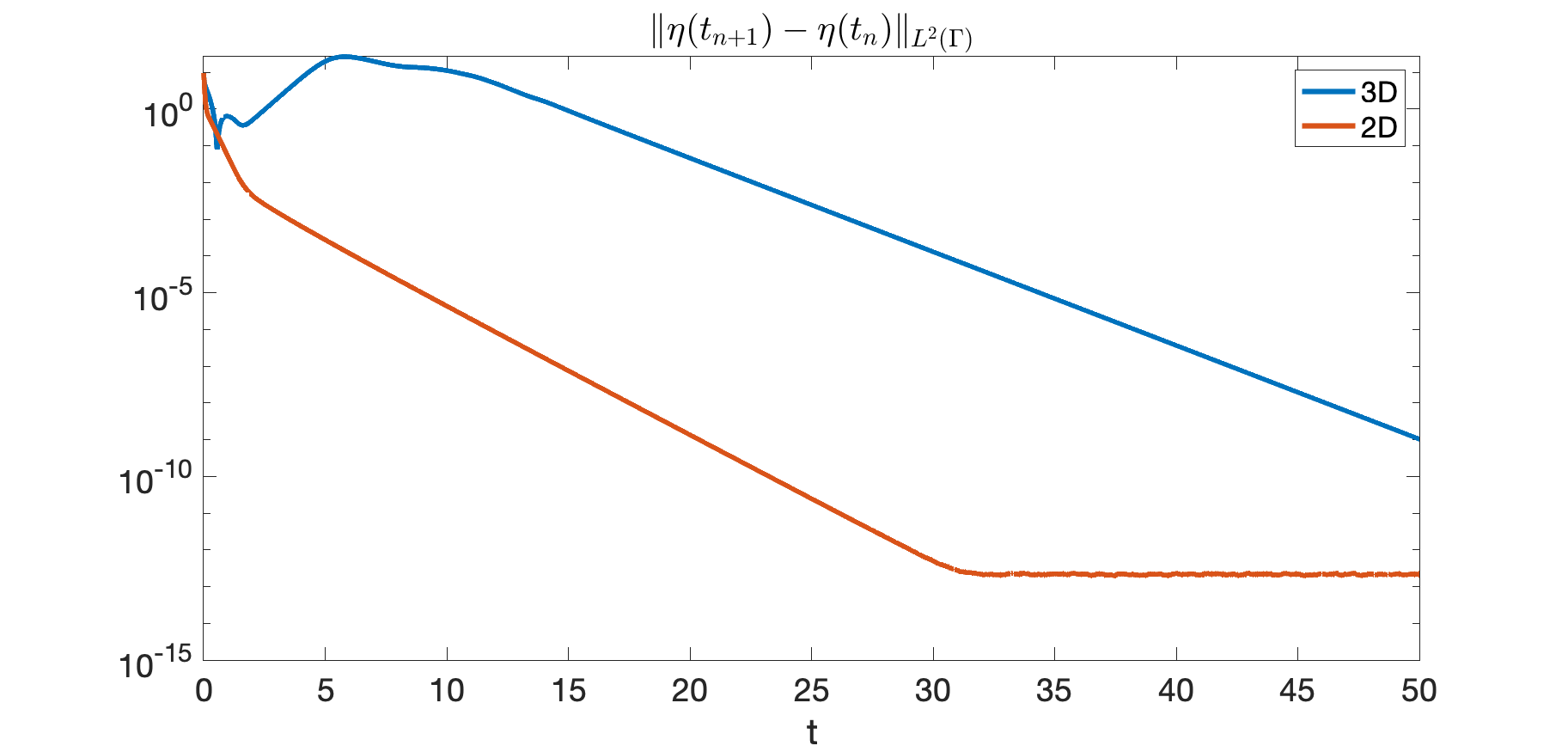}\label{fig:experimentT4_increment}}
\caption{Simulation T4.  Bifurcation and coupling  parameters as in Table \ref{tab:experiments_recap}. The BSDIB model \eqref{model} and the DIB model \eqref{model2d} reach different homogeneous steady states.}
\label{fig:experimentT4}
\end{center}
\end{figure}

\begin{figure}[h]
\begin{center}
\subfigure[Surface component $\eta$ at the final time $T=200$. ]{\includegraphics[scale=.4]{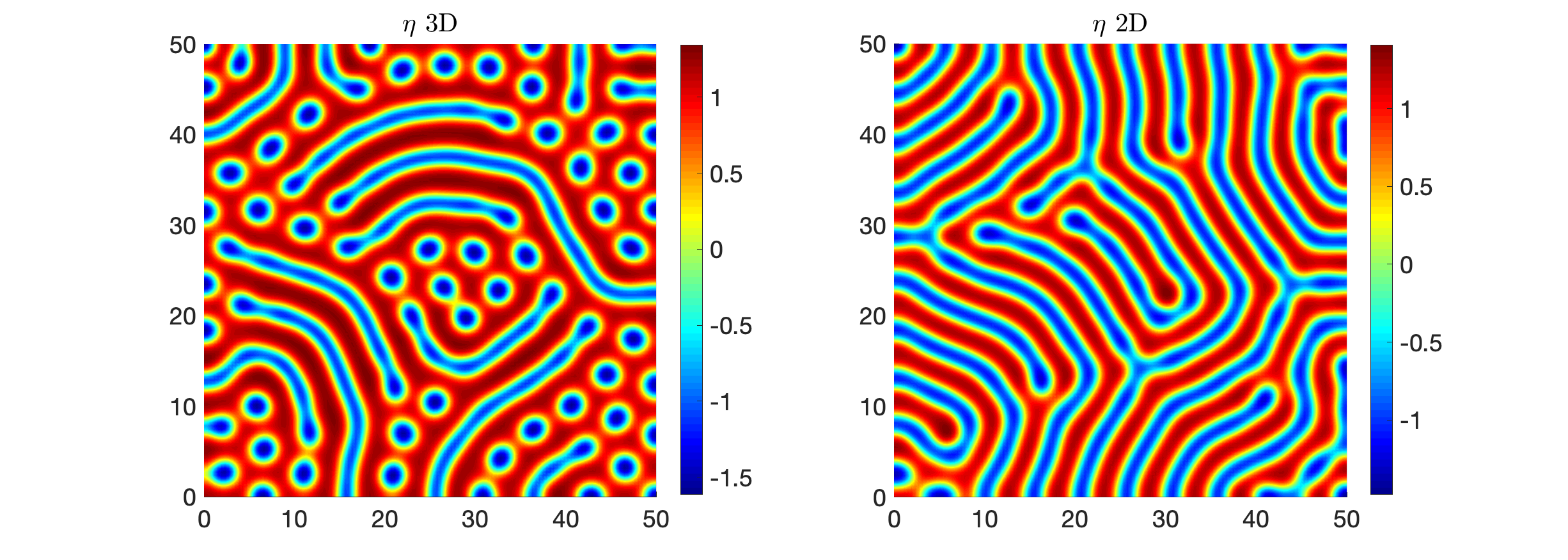}\label{fig:experiment1bis_comparison_T100}}
\subfigure[Surface component $\theta$ at the final time $T=200$. ]{\includegraphics[scale=.4]{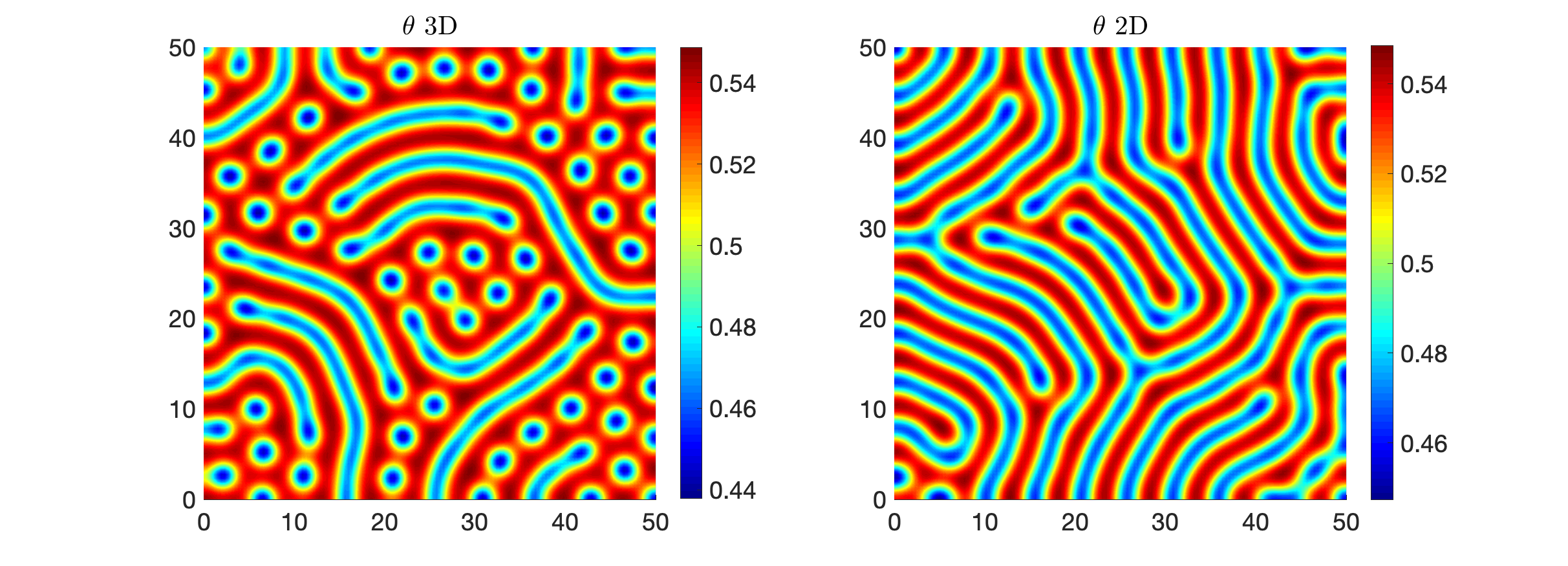}\label{fig:experiment1bis_comparison_T200}}
\subfigure[BS-DIB model \eqref{model}: bulk component $b$ and surface component $\eta$. ]{\includegraphics[scale=.4]{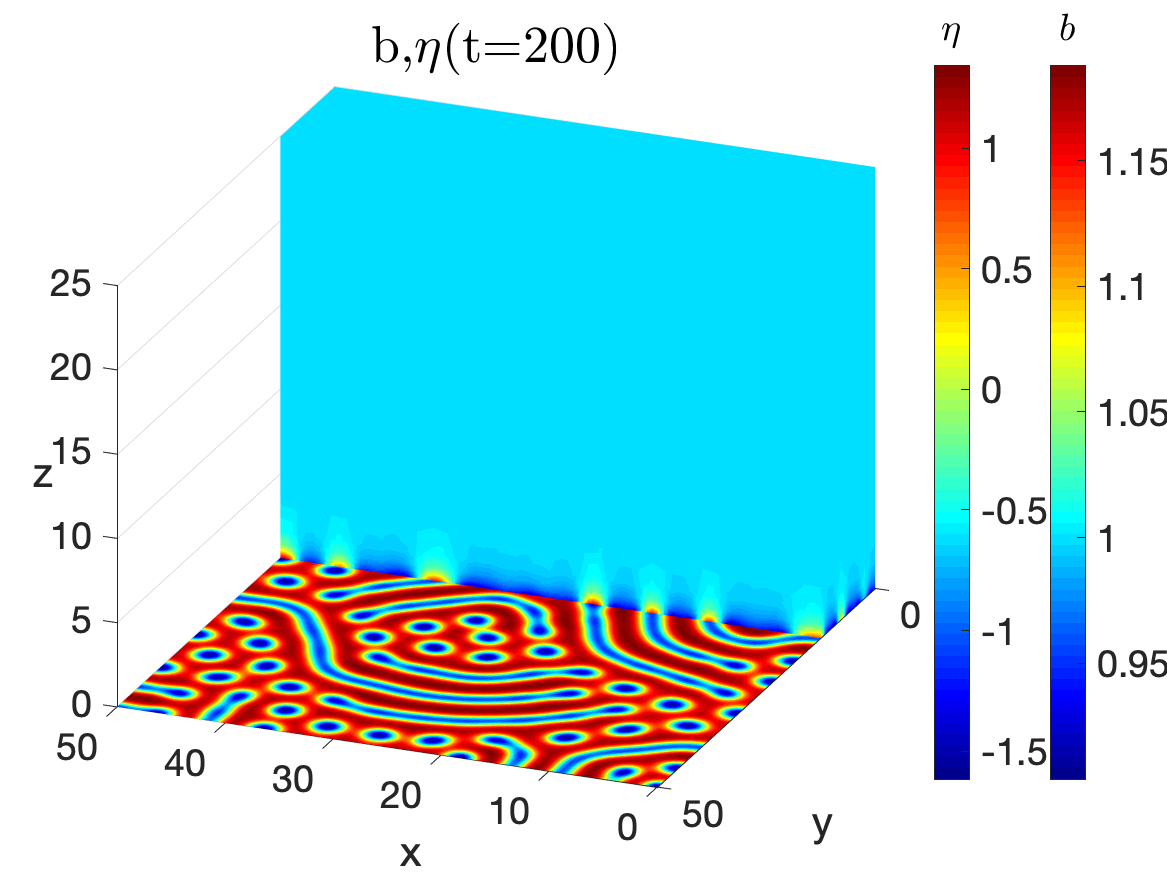}\label{fig:experiment1bis_bs}}
\subfigure[Increment of surface component $\eta$ over time.]{\includegraphics[scale=.4]{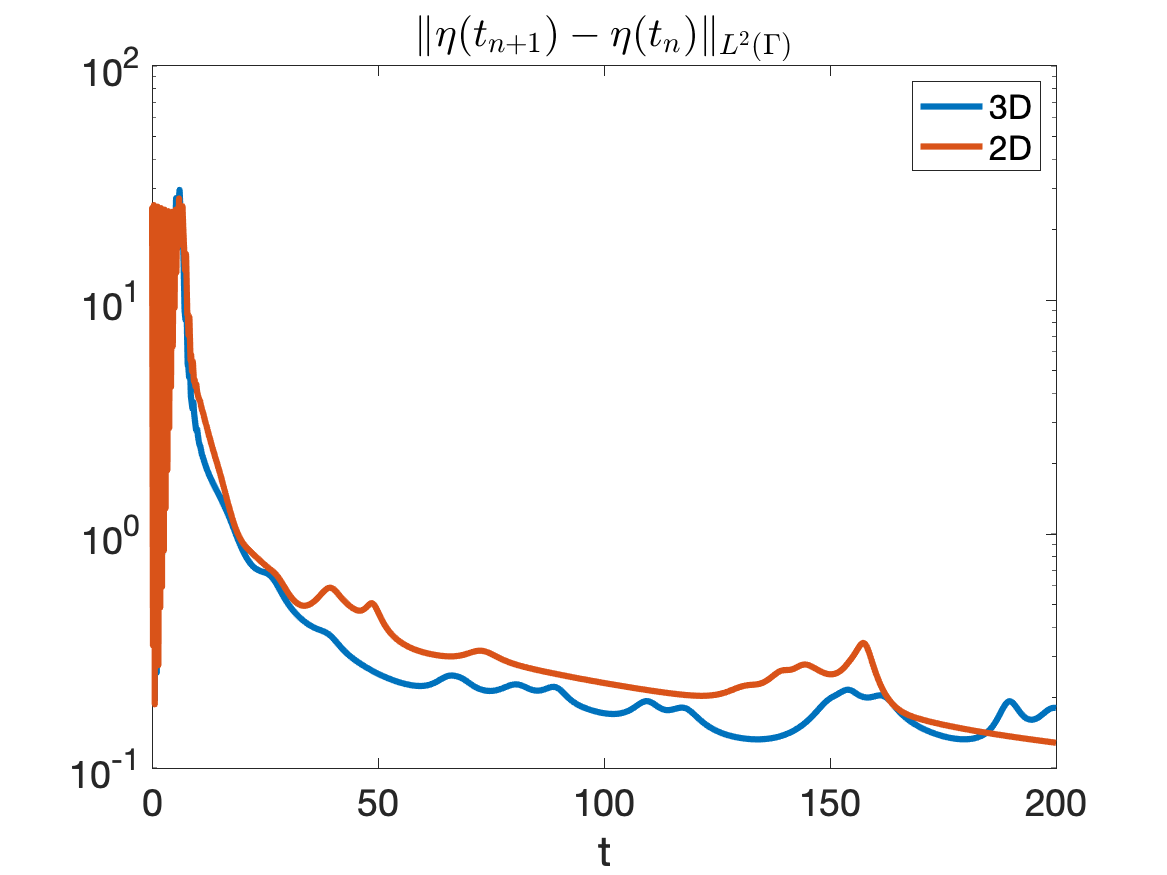}\label{fig:experiment1bis_increment}}
\caption{Simulation D1.  Bifurcation and coupling  parameters as in Table \ref{tab:experiments_recap}. The asymptotic steady state for the $\eta$ component shows a labyrinth pattern for the 2D DIB model \eqref{model2d} and a spot-and-worm pattern for the 3D BS-DIB model \eqref{model}. The spots of the spot-and-worm pattern in the coupled model \eqref{model} tend to slowly merge into worms over time.}
\label{fig:experiment1bis}
\end{center}
\end{figure}

\begin{figure}[h]
\begin{center}
\subfigure[Surface component $\eta$ at the final time $T=50$. ]{\includegraphics[scale=.4]{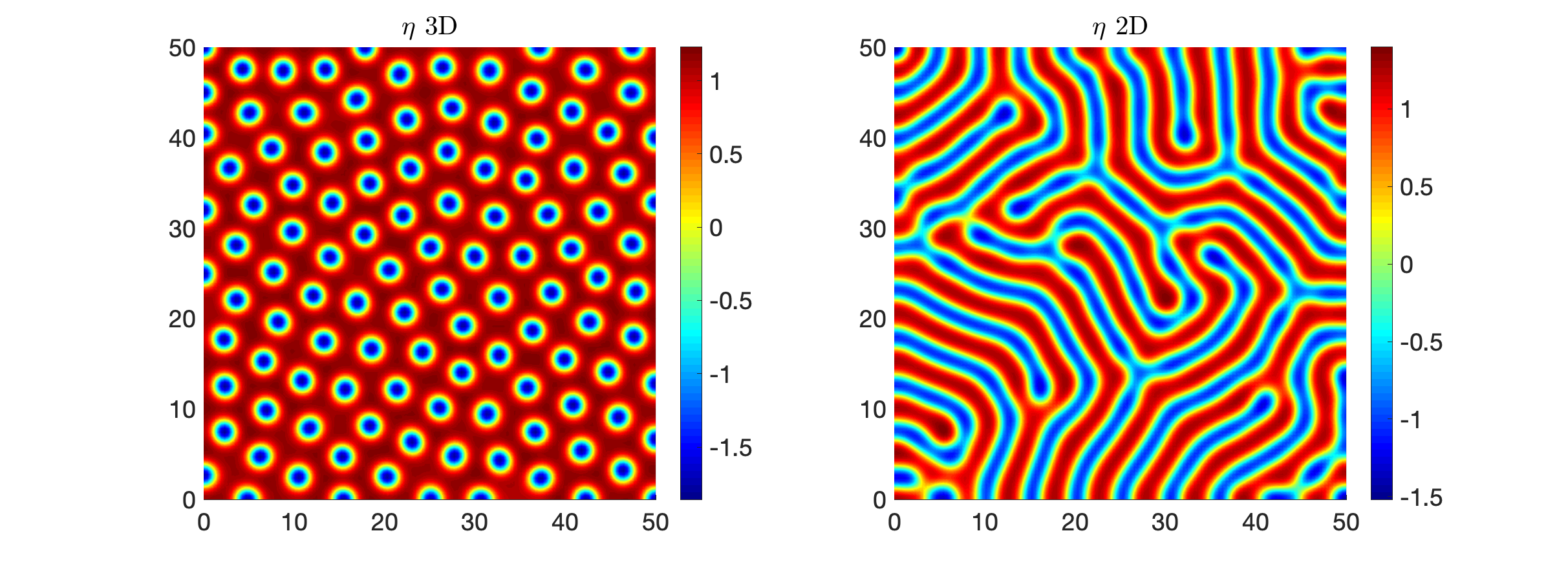}\label{fig:experiment1bis_comparison_psi02_eta}}
\subfigure[Surface component $\theta$ at the final time $T=50$.]{\includegraphics[scale=.4]{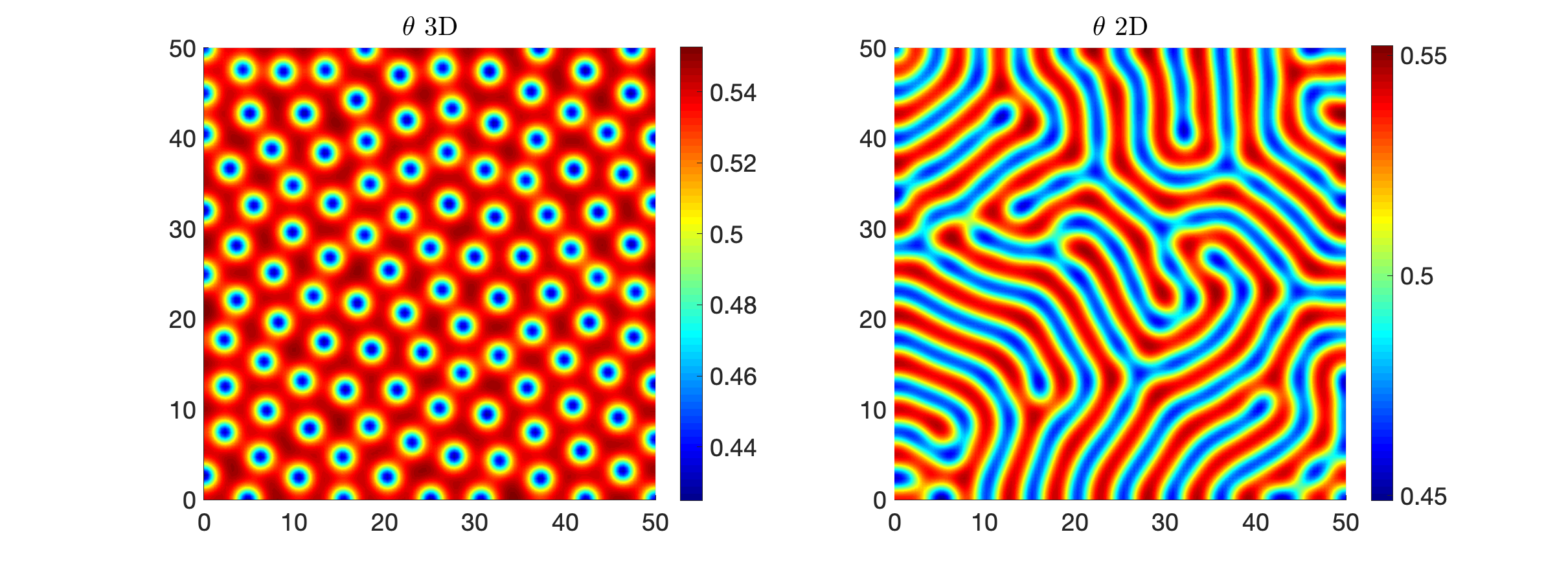}\label{fig:experiment1bis_comparison_psi02_theta}}
\subfigure[Increment of surface component $\eta$ over time.]{\includegraphics[scale=.4]{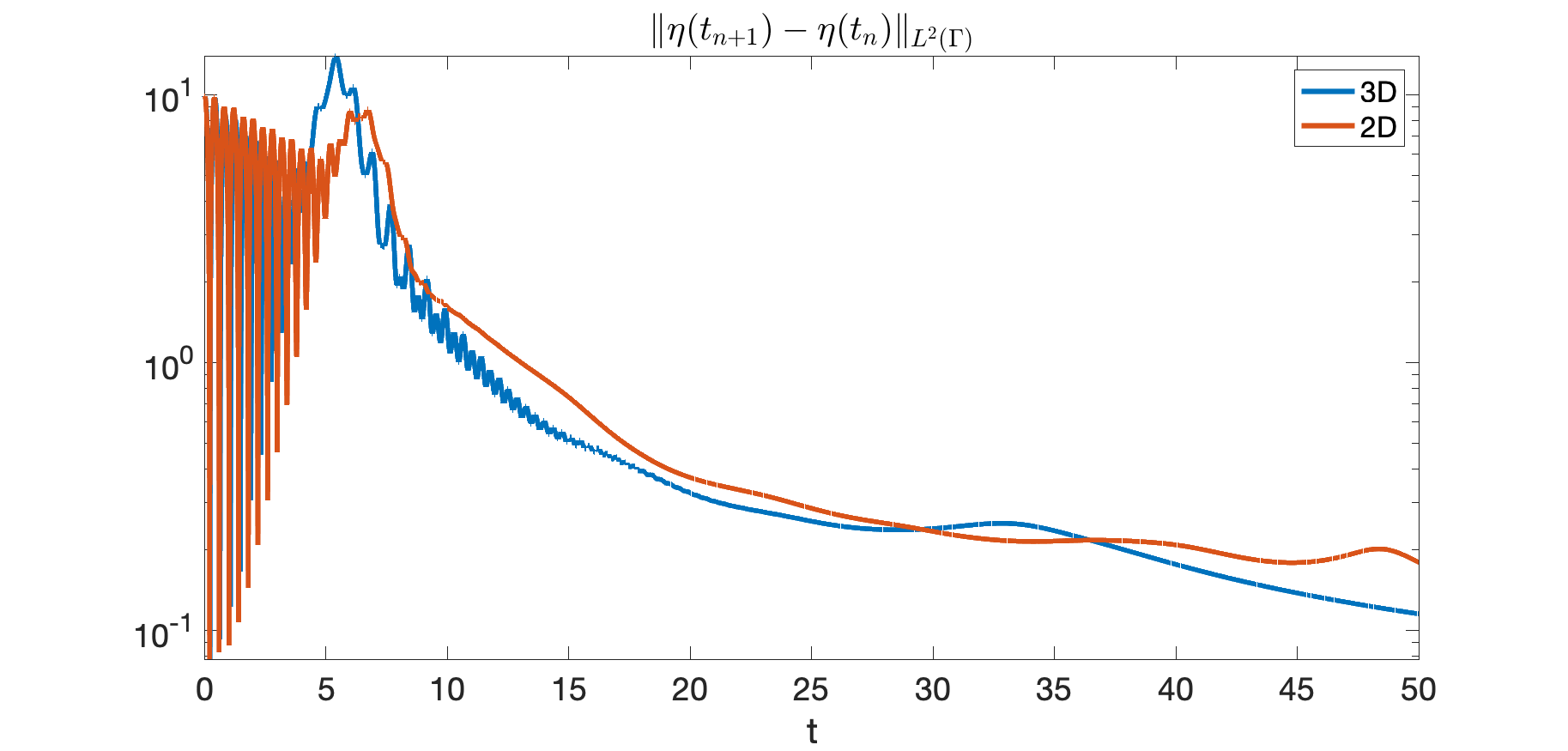}\label{fig:experiment1bis_psi02_increment}}
\caption{Simulation D3.  Bifurcation and coupling  parameters as in Table \ref{tab:experiments_recap}.  The 2D DIB model \eqref{model2d} exhibits a labyrinth pattern, while the 3D BS-DIB model \eqref{model} exhibits a reversed spots pattern. The results are in good agreement with those obtained with the MO-FEM approach in Fig. \ref{fig:experimentD3_MOFEM}, since both methods produce a spatial pattern of the same morphological class.}
\label{fig:experiment1bis_psi02}
\end{center}
\end{figure}

\begin{figure}[h]
\begin{center}
\subfigure[Surface component $\eta$ at the final time $T=200$. ]{\includegraphics[scale=.4]{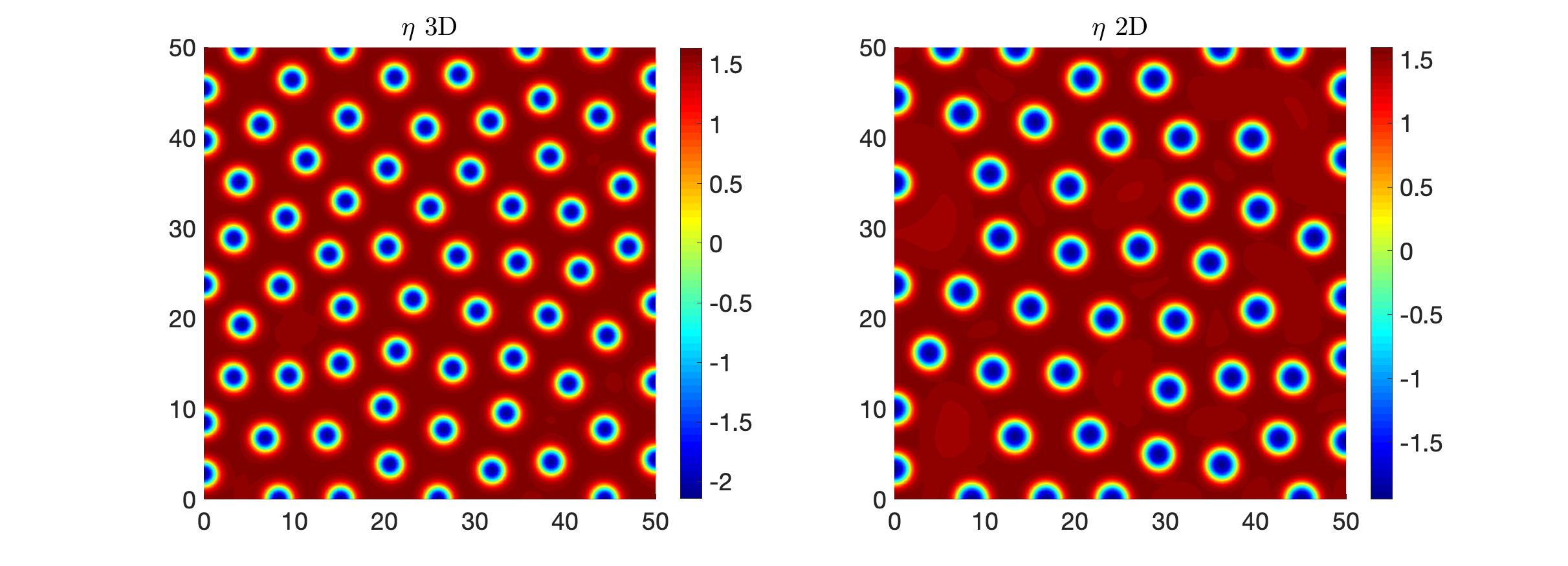}\label{fig:experiment2_comparison_eta}}
\subfigure[Surface component $\theta$ at the final time $T=200$.]{\includegraphics[scale=.4]{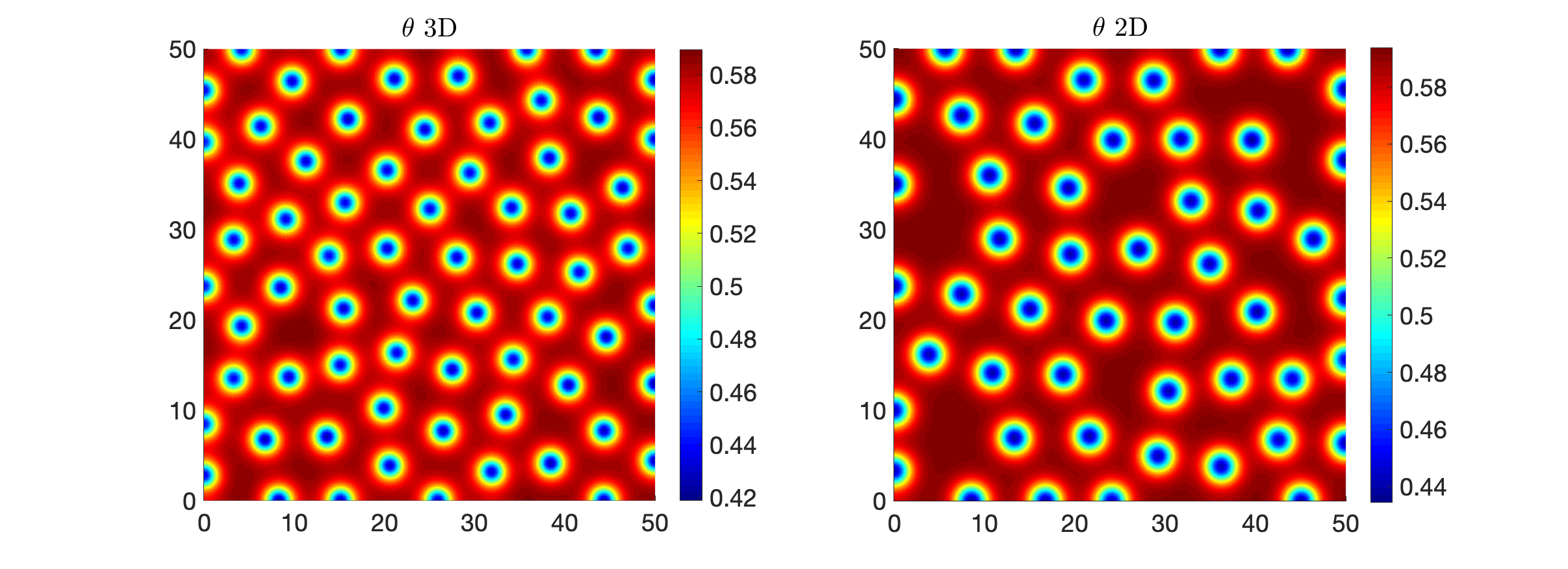}\label{fig:experiment2_comparison_theta}}
\subfigure[Increment of surface component $\eta$ over time.]{\includegraphics[scale=.4]{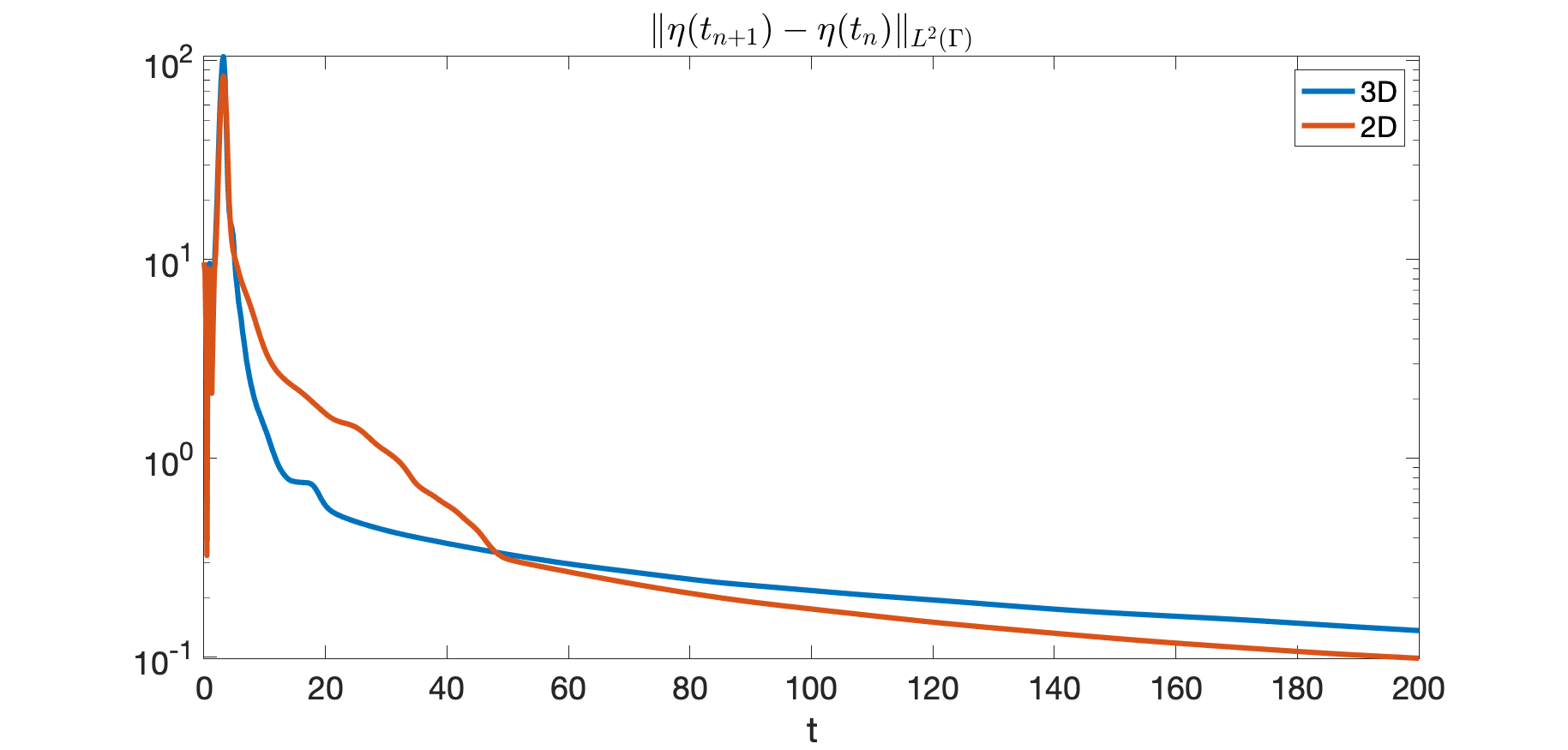}\label{fig:experiment2_increment}}
\caption{Simulation D4.  Bifurcation and coupling  parameters as in Table \ref{tab:experiments_recap}. The asymptotic steady state for the $\eta$ component shows reversed spots in both models. However, the coupled 3D BS-DIB model \eqref{model} shows smaller reversed spots with higher spatial density than the 2d BIB model \eqref{model2d}.}
\label{fig:experiment2}
\end{center}
\end{figure}

\begin{figure}[h]
\hspace*{-19mm}
\includegraphics[scale=0.39]{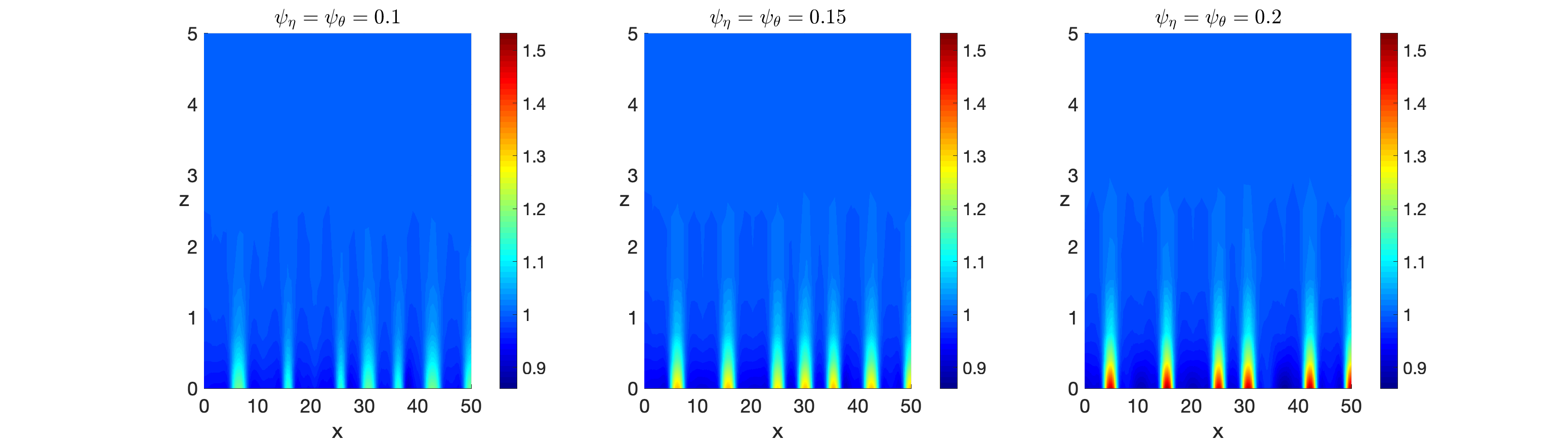}
\caption{Comparison of a section of the bulk component $b$ of Experiments D1, D2, D3 (see Table \ref{tab:experiments_recap}) for $y=12.5$. By increasing the coupling parameters $\psi_\eta = \psi_\theta$ (left to right), the bulk patterns increase in range (as indicated by the colorbars) and length along the $z$-direction.}
\end{figure}

\FloatBarrier

\section{Conclusions}
\label{sec:conclusions}
We have introduced a BSRD model in 3D, which we have called BS-DIB model, for electrodeposition. Compared to the previous DIB model in 2D, the new model fully accounts for the non-uniform electrolyte concentration in a neighborhood of the electrodic surface. The two-way coupling between bulk and surface substantially influences the long-term behavior of the system and in particular the morphological class of the Turing patterns obtained as asymptotic steady state solutions. Specifically,  we find that the bulk-surface coupling has two main effects.  First, we observe empirically that the BS-DIB model possesses a large Turing region in the parameter space, compared to the DIB model. Second,  when the parameters are chosen in the Turing space of the DIB model,  the BS-DIB model still exhibits spatial patterns, but of a different morphological class, i.e. the bulk-surface coupling affects the morphological class of the attained patterns.

The BS-DIB model is posed on a cubic domain, so it lends itself to efficient numerical solvers specifically devised for Cartesian grids, such as the MO-FEM.  Moreover, since the BS-DIB model exhibits spatial patterns only in a neighborhood of the surface,  we have adopted the BS-VEM on a graded mesh that is highly refined close to the surface and much coarser away from the surface.  Such a graded mesh combines the advantages of (i) being composed by equal elements of cubic shape, which significantly speeds up matrix assembly and improves matrix structure and (ii) has far less degrees of freedom than a uniform Cartesian grid with the same level of refinement close to the surface.  For this reason, the BS-VEM on a graded mesh proves to be more computationally efficient than the MO-FEM and is thus the spatial method of choice throughout this work.

As opposed to the MO-FEM, which is confined to structured geometries such as Cartesian grids, the BS-VEM can handle domains of general shape, thereby facilitating the simulation of real case studies.
A theoretical Turing instability analysis of the BS-DIB model is beyond the scope of this work. These aspects form part of our current investigations.






\section*{Acknowledgments and funding}
The work of MF was funded by Regione Puglia (Italy) through the research programme REFIN-Research for Innovation (protocol code 901D2CAA, project number UNISAL026) {and by the research project ``Metodi numerici innovativi per lo studio delle batterie'' (INdAM-GNCS project CUP\_ E55F22000270001)}.
The work of IS is supported by the MIUR Project PRIN 2020, ``Mathematics for Industry 4.0'', Project No. 2020F3NCPX, by the ICSC – Centro Nazionale di Ricerca in High Performance Computing, Big Data and Quantum Computing, funded by European Union – NextGenerationEU, Project code CN00000013 and by the research project ``Tecniche avanzate per problemi evolutivi: discretizzazione, algebra lineare numerica, ottimizzazione'' (INdAM-GNCS project CUP\_ E55F22000270001).
MF and IS are members of the INdAM-GNCS activity group (Italian National Group of Scientific Computing). MF is a member of SIMAI.

\section*{Conflict of interest}
The authors  declare no conflict of interest.

\FloatBarrier

\end{document}